\documentclass[12pt]{amsart} 
\usepackage[left=3cm,top=2.5cm,bottom=2.5cm,right=3cm]{geometry}
\usepackage{amssymb, amsmath, amsthm}
\usepackage{mathtools}
\usepackage{graphicx,xspace}
\usepackage{epsfig}
\usepackage{enumitem}
\usepackage[usenames,dvipsnames]{xcolor}
\usepackage{tikz}
\usepackage{stmaryrd}
\usepackage[T1]{fontenc}
\usepackage[utf8]{inputenc}
\usepackage{bm}
\usepackage[config, labelfont={normalsize}]{caption,subfig}

\usepackage{mathrsfs}
\definecolor{green}{RGB}{0,127,0}
\definecolor{red}{RGB}{191,0,0}
\usepackage[colorlinks,cite color=red,link color=green,pagebackref=true]{hyperref}

\usepackage{todonotes}

\makeatletter
\renewcommand{\tocsection}[3]{%
	\indentlabel{\@ifnotempty{#2}{\bfseries\ignorespaces#1 #2.\quad}}\bfseries#3}
\renewcommand{\tocsubsection}[3]{%
	\hspace{1cm} \indentlabel{\@ifnotempty{#2}{\ignorespaces#1 #2.\quad}}#3}
\makeatother

\usepackage[capitalize]{cleveref}

\theoremstyle{plain}
\newtheorem{lemma}{Lemma}
\newtheorem{theorem}{Theorem}
\newtheorem{corollary}{Corollary}
\newtheorem{proposition}{Proposition}

\newtheorem{definition}{Definition}

\theoremstyle{remark}

\newcommand{\QQ}{\mathbb{Q}}

\newcommand{\bM}{\bm{M}}
\newcommand{\bN}{\bm{N}}

\newcommand{\yy}{\bm{y}}
\newcommand{\pp}{\bm{p}}
\newcommand{\qq}{\bm{q}}
\newcommand{\rr}{\bm{r}}

\DeclareMathOperator{\sgn}{sgn}

\title[Constraints for $b$-deformed constellations]{Constraints for $b$-deformed constellations}
\author[V.~Bonzom]{Valentin Bonzom}
\address{LIGM, CNRS UMR 8049, Université Gustave Eiffel, Champs-sur-Marne, France}
\email{bonzom@univ-eiffel.fr}

\author[V.~Nador]{Victor Nador}
\address{LaBRI, UMR CNRS 5800, Universit\'e de Bordeaux, 351 cours de la Lib\'eration, 33405 Talence, France, EU \& Universit\'e Sorbonne Paris Nord, LIPN, CNRS, UMR 7030, F-93430 Villetaneuse, France}
\email{victor.nador@u-bordeaux.fr}

\thanks{VB is partially supported by the ANR-20-CE48-0018 3DMaps and the ANR-21-CE48-0017 LambdaComb. VN is supported by {\it Narodowe Centrum Nauki}, grant 2021/42/E/ST1/00162. In addition, the authors are thankful to Guillaume Chapuy and Maciej Do\l\k{e}ga for numerous explanations of their main work which our article is based on.}

\begin{document}
\maketitle

\begin{abstract}
Hurwitz numbers count branched covers of the sphere and have been of interest in various fields of mathematics. Motivated by the Matching-Jack conjecture of Goulden and Jackson, Chapuy and Do\l\k{e}ga recently introduced a notion of $b$-deformed double weighted Hurwitz numbers. It equips orientable and non-orientable maps and constellations with $b$-weights defined inductively. It is then unclear whether some elementary properties of orientable maps remain true due to the nature of the $b$-weights. We consider here the case of the Virasoro constraints, which express that choosing a corner is equivalent to rooting in terms of generating functions. We prove this property for two families of $b$-deformed Hurwitz numbers, namely 3-constellations and bipartite maps with black vertices of degrees bounded by 3. The proof is built upon functional equations from Chapuy and Do\l\k{e}ga and a lemma which extracts the constraints provided they close in some appropriate sense for the commutator. This requires to calculate the constraint algebra which in those two families do not form a Lie algebra but a generalization of independent interest with structure operators instead of structure constants.
\end{abstract}

\section*{Introduction}
\subsection*{$b$-deformed Hurwitz numbers and the Matching-Jack conjecture}
Hurwitz numbers \cite{Hurwitz1891} count the numbers of inequivalent branched covers of the 2-sphere by surfaces of arbitrary genus with some given ramification profiles. In particular, \emph{double Hurwitz numbers} count coverings with two fixed profiles given by some partitions $\lambda, \mu$, say over 0 and $\infty$, and $k$ other ramification points whose profiles are not recorded but whose number of preimages is tracked instead. The genus is then given by the Riemann-Hurwitz formula. Finding those numbers can be reformulated as a purely combinatorial problems in terms of transitive factorizations in the symmetric group \cite{Frobenius1900}. This reveals a deep connection to algebraic combinatorics which has been thoroughly investigated with impressive results. In particular, the generating function of (connected or not) double Hurwitz numbers can be expressed as a sum of products of Schur functions $s_\lambda$,
\begin{equation}
\tau_k(t, \pp, \qq, u_1, \dotsc, u_k) =  \sum\limits_{n \geq 0} t^n \sum_{\lambda \vdash n} s_\lambda(\pp) s_\lambda(\qq) \prod_{\square\in\lambda} \prod_{i=1}^k (u_i + c(\square)),
\end{equation}
which shows directly \cite{HarnadOrlov2015, GuayPaquetHarnad2015, Guay-PaquetHarnad2017} that it is a tau function of the Toda hierarchy, of hypergeometric type as introduced by Orlov and Scherbin \cite{OrlovScherbin2000}.

The above generating function is an instance of \emph{weighted} Hurwitz numbers \cite{HarnadOrlov2015, GuayPaquetHarnad2015, Guay-PaquetHarnad2017}. It is said to have polynomial weight function $G(z) = \prod_{i=1}^k (u_i + z)$, which in the case of weighted Hurwitz numbers can be replaced with a function in a more general class. Weighted Hurwitz numbers have attracted a lot of attention, due to their numerous interesting combinatorial properties. They range from closed formulas for the enumeration of planar coverings \cite{BousquetSchaeffer2000}, to the topological recursion \cite{AlexandrovChapuyEynardHarnad2020, BychkovDuninBarkowskiKazarianShadrin2020, BCCGF22, BDBKS4} for correlators with arbitrary genus and marked faces.

The above expression in terms of Schur functions has a natural generalization from the viewpoint of algebraic combinatorics, by deforming the Schur functions into Jack symmetric functions $J^{(\alpha)}_\lambda$ and the content of a box $c(\square)$ into its $\alpha$-deformed content \cite{Jack1970/1971, Stanley1989},
\begin{equation} \label{GouldenJacksonSeries}
\tau^{(\alpha)}_k(t, \pp, \qq, u_1, \dotsc, u_k) =  \sum\limits_{n \geq 0} t^n \sum_{\lambda \vdash n} \frac{J^{(\alpha)}_\lambda(\pp) J^{(\alpha)}_\lambda(\qq)}{j_\lambda^{(\alpha)}} \prod_{\square\in\lambda} \prod_{i=1}^k (u_i + c_\alpha(\square)).
\end{equation}
Studying this deformation was a starting point of Goulden and Jackson's work in \cite{GouldenJackson1996Jack}. They considered not double Hurwitz numbers like above, but the so-called deformed hypermaps series $\tau^{(\alpha)}(\pp, \qq, \rr) = \sum_\lambda \frac{1}{j_\lambda^{(\alpha)}} J^{(\alpha)}_\lambda(\pp) J^{(\alpha)}_\lambda(\qq) J^{(\alpha)}_\lambda(\rr)$, which corresponds to Hurwitz numbers with exactly three ramification points, whose profiles are fixed partitions $\lambda, \mu, \nu$. They remarkably found (but did not prove) that in the basis $p_\lambda q_\mu r_\nu$, the coefficients are \emph{polynomials in $b=\alpha-1$ with non-negative integer coefficients} (up to a known denominator). This is presently a conjecture, which they called the Matching-Jack conjecture. To this day, polynomiality has been proven in~\cite{DolegaFeray2016} and integrality in~\cite{BenDali:Integrality}, but the positivity of the coefficients is still missing.

An analogous conjecture can be formulated for the series $\tau^{(\alpha)}_k(t, \pp, \qq, u_1, \dotsc, u_k)$, which was proved to be true by Chapuy and Do\l\k{e}ga in \cite{ChapuyDolega2020}, i.e. the coefficients of $\tau^{(\alpha)}_k(t, \pp, \qq, u_1, \dotsc, u_k)$ in the basis $p_\lambda q_\mu \prod_{i=1}^k u_i^{l_i}$ are polynomials in $b$ with integer coefficients (up to a known denominator). The combinatorial idea behind \cite{ChapuyDolega2020} lies in the interpretation of Hurwitz numbers as \emph{constellations}, a generalization of combinatorial maps. While this is a well-known correspondence for $b=0$, laid down by Lando and Zvonkin \cite{LandoZvonkin2004}, it was lacking for other values of $b$.

\subsection*{Non-oriented constellations}
Due to cases of \eqref{GouldenJacksonSeries} at $b=0$ and $b=1$ having been previously understood as counting respectively orientable maps and constellations \cite{LandoZvonkin2004, Guay-PaquetHarnad2017, HarnadOrlov2015} and non-oriented\footnote{Throughout this article, \emph{non-oriented} means ``orientable or not''.} maps \cite{GouldenJackson1996Zonal,Dali21}, it is natural to expect the parameter $b$ to be a sort of measure of non-orientability on discrete surfaces. This idea was put forward by La Croix \cite{LaCroix2009} who used it to prove the most general case of the $b$-conjecture,
on ordinary maps counted by face degrees and number of vertices\footnote{The specialization of the hypermap series is: arbitrary $\pp$, but $\qq = (0, 1, 0, \dotsc)$ and $\rr = (u, u, \dotsc)$. Equivalently, an arbitrary ramification profile over $0$, only transpositions over $\infty$ and above the third point, one only keeps track of the number of preimages and not the full profile. Notice that it is also a special case of $\tau^{(\alpha)}_k(t, \pp, \qq, u_1, \dotsc, u_k)$ with $k=1$ and $\qq = (0, 1, 0, \dotsc)$.}.

In \cite{ChapuyDolega2020}, Chapuy and Do\l\k{e}ga introduced non-oriented weighted Hurwitz numbers. In the case of the polynomial function $G(z) = \prod_{i=1}^k (u_i + z)$, they correspond to \emph{non-oriented constellations} equipped with a \emph{measure of non-orientability}. It assigns a $b$-weight to every rooted, non-oriented constellations. That $b$-weight is defined inductively and it is not apparently symmetric under the exchange of the variables $\pp$ and $\qq$. This makes proving properties which are obvious for orientable maps and constellations far from obvious in the deformed setting, if they are true at all.

We prove in the present article one such property in two classes of models from~\cite{ChapuyDolega2020} defined below, namely that for these models, marking a corner has the effect of rooting at the level of the generating functions. This property was conjectured in~\cite[Eq.$(17)$]{ChapuyDolega2020} for general weighted Hurwitz numbers, and proved for monotone Hurwitz numbers in \cite{BonzomChapuyDolega2021}\footnote{In addition to general maps and bipartite maps where it was indirectly known by combining previous results relating the series \eqref{GouldenJacksonSeries} to the $\beta$-ensemble of matrix integrals \cite{Okounkov1997, AdlervanMoerbeke2001, LaCroix2009}.} for which the weight function is $G(z) = \frac{1}{u+z}$. A simple corollary is that the generating series and the one calculated with the $b$-weights of the dual constellations are equal, which is a non-trivial result as dicussed in~\cite{LaCroix2009}.

The technique we use here relies on proving that the operators involved by the rooting form a sort of algebra for the commutator (defined precisely below). In the case of ordinary maps, it is known to be a (half-)Virasoro algebra~\cite{AdlervanMoerbeke2001}, hence the name \emph{Virasoro constraints} for the equations they produce on the tau function. Remarkably, it is also a (half-)Virasoro algebra for $b$-deformed monotone Hurwitz numbers \cite{BonzomChapuyDolega2021}. But the extensions we find here are new algebras and we believe that they are of independent interest. The two classes of models we consider are generalizations of bipartite maps and general maps respectively: 
\begin{enumerate}[label=\Roman*]
\item \label{item:A}-- constellations with at most 3 colors: $\qq = (1,0, 0, \dotsc)$ and $k\leq 3$
\item \label{item:B}-- bipartite maps with controlled degrees bounded by 3 for black vertices: $\qq = (q_1,q_2,q_3,0, \dotsc)$ and $k=1$.
\end{enumerate}
In both cases, we obtain \emph{higher degree constraints} which are cubic with respect to the variables $\pp$ instead of quadratic as in the case of maps and bipartite maps. The resulting constraint algebras depend on $b$ explicitly and are not Lie algebras (they have ``structure operators'' instead of structure constants). A specialization to $b=0$ has already appeared in~\cite{MarshakovMironovMorozov}.

\subsection*{Obstacles}
Traditional approaches to the cases $b=0, 1$ are representation theory and matrix integrals. Indeed, at $b=0, 1$, orientable and non-oriented constellations can be encoded as factorizations in the symmetric group and in a double coset respectively, which via some Frobenius' formulas give rise to the Schur and zonal expansions directly \cite{GouldenJackson1996Zonal,Dali21}. However, there is no such representation theory to be used for general $b$, which is of course one of the interests in the matching-Jack conjecture. 

There is a long history of connections between maps and matrix integrals, but our results go beyond what is known to be achievable in terms of matrix integrals. In addition to the well-established connection between ordinary maps (respectively bipartite maps) and the Gaussian Unitary Ensemble (GUE) (respectively Laguerre ensemble), Okounkov \cite{Okounkov1997} related a specialization of the hypermap series to the Gaussian $\beta$-ensemble, which is a 1-parameter generalization of the GUE with $\beta = 2/(1+b)$, and this was the specialization later considered by La Croix for which he could prove the Matching-Jack conjecture. On the side of (orientable) constellations, there does exist a matrix model, by Ambj\o{}rn and Chekhov~\cite{AmbjornChekhov2014}, for orientable constellations. However, it is a multi-matrix model and as such it is not possible to use the standard approaches to matrix integrals based on reducing the integral to eigenvalues. 

Furthermore, it is not known how to deform multi-matrix models such as \cite{AmbjornChekhov2014}. Defining 2-matrix models for the $\beta$-ensemble has been an issue in itself \cite{BergereEynard2009, BergereEynardMarchalPrats-Ferrer2012}. One approach is actually to \emph{define} them in terms of Jack expansions (and not proving it from an integral definition) \cite{MironovMorozovShakirov2010}. An interesting question arising from our work is therefore whether there are matrix integrals for the models we present (see \cite{GisonniGravaRuzza2021,ruzza2023jacobi} for a promising approach using limits of the Jacobi $\beta$-ensemble).

Our approach instead relies on $(i)$ results from~\cite{ChapuyDolega2020} of combinatorial and algebraic origin, $(ii)$ a critical lemma from~\cite{BonzomChapuyDolega2021}, and $(iii)$ some quite heavy calculations to express the commutators between the higher order constraints of our models. Ingredients $(i)$ and $(ii)$ are fairly generic and hold for general $b$-deformed weighted Hurwitz numbers. Ingredient $(iii)$ is the limiting factor as the calculations are tricky and we have not found a way to bypass them or put them on a computer.

Below we describe the plan, setup and results more precisely.

\subsection*{Presentation of the results}
A map is a graph equipped with a cyclic ordering of the edges meeting at every vertex (loops and multiple edges are allowed). Let $k\geq 1$ be an integer, then orientable $k$-constellations are generalizations of bipartite maps to $k+1\geq 2$ colors. Edges are generalized to \emph{hyperedges}, i.e. $k+1$-gons whose vertices carry the colors from $0$ to $k$ in clockwise (or counter-clockwise) order around each $k+1$-gon. The size of a constellation is its number of hyperedges. Orientable $k$-constellations count branched covers to the 2-sphere with $k+2$ ramification points \cite{LandoZvonkin2004}.

In \cite{ChapuyDolega2020}, Chapuy and Do\l\k{e}ga introduced not only non-oriented constellations but also \emph{$b$-deformed constellations}\footnote{This is one of the first steps in their article, which corresponds to polynomial weight functions $G$. Then they extend their setup to weighted Hurwitz numbers. The combinatorial objects are then not called constellations anymore although their properties can be derived from those of constellations.}, which are non-oriented constellations equipped with a $b$-weight, where $b$ is a formal parameter. 
Each connected, rooted constellation $(\bM,c)$ ($c$ being the root, i.e. an oriented corner) receives a $b$-weight $\rho_b(\bM,c)$ which is a polynomial in $b$ and is equal to 1 for orientable constellations. It is evaluated inductively, through a deletion algorithm which starts from the root of $(\bM, c)$. If $\bM$ is a labelled, non-rooted constellation, its weight is defined as
\begin{equation} \label{AverageWeight}
\rho_b(\bM) = \frac{1}{n} \sum_{\text{rooting $c$ of $\bM$}} \rho_b(\bM,c)
\end{equation}
i.e. it is the average of the weights over all choices of root. For a non-connected constellation $\bM$, its weight is defined multiplicatively over its connected components.

Deformed constellations still have a notion of faces and their topologies are in fact those of non-oriented surfaces. Chapuy and Do\l\k{e}ga introduced a generating function which controls the degrees of faces, as well as the degrees of vertices of color 0, and the number of vertices of every other colors. Denote $\tau$ this series for non-necessarily connected constellations, with parameters $t$ tracking the size, $p_i$ for faces of degree $i$, $q_j$ for vertices of color 0 and degree $j$. Let also $H = (1+b) \ln\tau$ by the generating series of connected constellations, and $H^{[m]}$ the generating function of connected, rooted constellations $(\bM,c)$ with a root vertex of degree $m$ and weighted with $\rho_b(\bM,c)$.

\begin{theorem} \cite{ChapuyDolega2020}
One has
\begin{equation} \label{VertexInvariance}
k\frac{\partial H}{\partial q_k} = H^{[k]}
\end{equation}
The generating function $\tau$ satisfies an evolution equation, i.e. an equation of the form
\begin{equation}\label{EvolutionEquation_CD}
\frac{\partial\tau}{\partial t} = D \tau
\end{equation}
where $D$ is a differential operator w.r.t. the variables $p_i$.
\end{theorem}

In~\cite{ChapuyDolega2020}, the evolution equation~\eqref{EvolutionEquation_CD} is derived from~\eqref{VertexInvariance} and heavy calculations involving Jack symmetric functions. Equation~\eqref{VertexInvariance} is interpreted as follows in~\cite{ChapuyDolega2020}. Consider the set of connected constellations $\bM$ which have at least one vertex of color 0 and degree $k$. Then the expectation of $\rho_b(\bM)$ over that set is the same as the expectation of $\rho_b(\bM,c)$ over connected, rooted constellations conditioned to have a root vertex of degree $k$. This equality is far from obvious due to the inductive nature of the $b$-weight calculated from the root.

\medskip

The two families \ref{item:A} and \ref{item:B} we consider are special cases of the $b$-deformed constellations of Chapuy and Do\l\k{e}ga \cite{ChapuyDolega2020}. Whenever we need to refer to the objects of our two families generically we will call them constellations. Our results will follow from the following theorem.
\begin{theorem}
For the two families of models \ref{item:A} and \ref{item:B}, the generating function $\tau$ is determined by a set of differential constraints, for $i\geq 1$
\begin{equation}
L_i \tau = 0, \quad \text{with} \quad L_i = i\frac{\partial}{\partial p_i} - t^r M_i,
\end{equation}
where $r>0$ is a positive integer and the $M_i$s are polynomials in $t$ whose coefficients are explicit differential operators in the variables $p_j$s. The constraints satisfy
\begin{equation} \label{ConstraintAlgebra}
[L_i, L_j] \coloneqq L_iL_j-L_jL_i = t\sum_{k \geq 1} D_{ij}^k L_{k}.
\end{equation}
The differential operators $D_{ij}^k$ are explicit for both families \ref{item:A} and \ref{item:B}, and given in Theorems \ref{thm:BipMaps} for bipartite maps (half-Virasoro), \ref{thm:3-constellations} for 3-constellations and \ref{thm:2ndModel} for bipartite maps with controlled black vertex degrees bounded by 3.
\end{theorem}

For ordinary maps 
and for bipartite maps 
, the $D_{ij}^k$s are constant with respect to the $p_i$s and the constraints form a half-Virasoro algebra (in particular, it is independent of $b$) as expected. 
In the orientable setting ($b=0$), a subfamily of \ref{item:B}, bipartite maps with black vertices having fixed degree equal to 3, was worked out in~\cite{MarshakovMironovMorozov} (although the argument is seemingly not mathematically correct). As far as we know, the commutators for 3-constellations had never appeared in the literature (the constraints themselves had originally been written by Fang~\cite{Fang:PhD} at $b=0$, before the much more general case of \cite{ChapuyDolega2020}).

In the orientable case there is an obvious symmetry between the face degrees and the degrees of vertices of color 0. It allows to directly interpret $k \frac{\partial H }{\partial p_k}$ as the generating series of maps with a rooted face of degree $k$. Here it is no longer immediate due to the way the $b$-weight is defined  (see~\cite[Thm. $3.35$]{LaCroix2009}). One consequence of our results is that despite the additional difficulties due to the nature of the $b$-weight, this interpretation holds in the two families of models we consider.

\begin{corollary} \label{thm:Rooting}
For the families \ref{item:A} and \ref{item:B}, denote $H^{(i)}$ the generating function of connected, rooted constellations $(\bM,c)$ with root face of degree $i$ and weighted with $\rho_b(\bM,c)$. Then
\begin{equation}
i\frac{\partial H}{\partial p_i} = H^{(i)}.
\end{equation}
\end{corollary}

Another corollary can be obtained by combining the above with a result of \cite{ChapuyDolega2020}. Let $\tilde{H}^{[i]}$ be the generating function of connected, rooted constellations $(\bM,c)$ with a root face of degree $i$, but for which the $b$-weight is $\rho_b(\tilde{\bM}, \tilde{c})$, calculated with the dual rooted constellation $(\tilde{\bM}, \tilde{c})$ to $(\bM, c)$.

\begin{corollary} \label{thm:Duality}
For the same models, $\tilde{H}^{[i]} = H^{(i)}$.
\end{corollary}
This formalizes and proves a conjecture of La~Croix \cite[Conjecture 3.33]{LaCroix2009} (formulated in the context of ordinary maps), motivated by the symmetry of $H$ under the exchange $q_i\leftrightarrow p_i$.

While we believe that all our results hold for any $b$-deformed, weighted Hurwitz numbers, our method, which partially relies on bruteforce computation of the constraint algebra, makes it too difficult to implement beyond present article.

\subsection*{Plan and notation}
Section \ref{sec:Lemma} presents a technical lemma extended one originally from \cite{BonzomChapuyDolega2021}. It generally enables one to extract differential constraints from ``evolution equation'' which is precisely what is needed to exploit the results of \cite{ChapuyDolega2020}. In Section \ref{sec:Definitions}, we recall everything we need from \cite{ChapuyDolega2020} and write the precise versions of our theorems and the short versions of the proofs which apply our lemma from Section \ref{sec:Lemma}. For it to be applicable, we need to check the constraint algebra of each model. We do so in one class of models in Section \ref{sec:3Constellations} (for 3-constellations), and in the second class of models in Section \ref{sec:commk1} (bipartite maps with vertices of color 0 having degrees bounded by 3).

In the rest of the article, we write $\pp \coloneqq (p_1, p_2, \ldots)$ an infinite set of indeterminates. 
If $t_1, \ldots, t_r$ are formal variables, and $R$ is a ring, we denote $R[[t_1, \ldots, t_r]]$ the ring of formal power series in  $t_1, \ldots, t_r$ with coefficients in $R$, e.g. $\QQ[\pp][[t]]$ is the ring of formal power series in $t$ whose coefficients are polynomials in the $p_i$s.

\section{The constraint lemma} \label{sec:Lemma}

We start with a slightly technical lemma, which extends Lemma $2.6$ from \cite{BonzomChapuyDolega2021}\footnote{The Lemma as stated in the published version is incorrect as it misses an assumption. The statement and its proof are corrected in the v3 updated on arXiv.}, and enables us to extract differential constraints from an ``evolution equation''. Here the variable $p_i$ has degree $i$ and the degree of a product of them is the sum of their degrees.

\begin{lemma}
\label{lemma:constr_from_ev}
    Let $Z\in \QQ[\pp][[t]]$ and let $r\geq 1$ be an integer. Assume that 
  \begin{enumerate}[label=(\roman*)]
  \item\label{enum:Homogeneity} $[t^n]Z$ is a homogeneous polynomial of degree $n$ in the variables $p_i$s.
   \item\label{enum:Evolution} there exist operators $M^{(1)}, \ldots, M^{(r)} \in \QQ[[p_1, p_2, \dotsc, \frac{\partial}{\partial p_1}, \frac{\partial}{\partial p_2}, \dotsc]]$ which are independent of $t$ such that
   \begin{equation} \label{EvolutionEquation}
   \frac{\partial Z}{\partial t} = \sum_{m=1}^r t^{m-1} M^{(m)} Z
   \end{equation}
   \item\label{enum:Algebra} For $m=1, \ldots, r$, $M^{(m)}$ can be rewritten as $M^{(m)} = \sum\limits_{i \geq 1} p_i M^{(m)}_i$, such that the operators
   \begin{equation}
   L_i \coloneqq i\frac{\partial}{\partial p_i} - \sum_{m=1}^r t^m M^{(m)}_i \ \in \QQ[[p_1, p_2, \dotsc, \frac{\partial}{\partial p_1}, \frac{\partial}{\partial p_2}, \dotsc]][t]
   \end{equation}
    for $i\geq 1$ satisfy the following two conditions
    	\begin{itemize}
    	\item $[t^n]L_i F$ is a homogeneous polynomial of degree $n-i$ in the variables $p_i$s,
    	\item the operators $L_i$ are closed under commutation in the sense that there exist $D_{ij}^k\in \QQ[[p_1, p_2, \dotsc, \frac{\partial}{\partial p_1}, \frac{\partial}{\partial p_2}, \dotsc]][t]$, for $i, j, k\geq 1$, such that
   		\begin{equation} \label{CommutationRelation}
   		\left[L_i,L_j\right] = t \sum_{k \geq 1} D_{ij}^k L_k
   		\end{equation}
   		\end{itemize}
 \end{enumerate}
    then $Z$ satisfies the constraints
    \begin{equation}
        \label{eq:constraints}
        L_i Z = 0 \qquad \text{for all $i\geq1$.}
    \end{equation}
\end{lemma}


Notice that the factor $t$ in the RHS of \eqref{CommutationRelation} is crucial. Without it, the commutation relation $\left[L_i,L_j\right] = \sum_{k \geq 1} D_{ij}^k L_k$ has the simple solution $D_{ij}^k = \delta_{j,k} L_i - \delta_{i,k} L_j$, which is obviously not interesting. As we will see, the difficulty in applying the above lemma is identifying the $M_i^{(m)}$s and actually proving the existence of the commutation relation \eqref{CommutationRelation}.

As another technical remark, notice that in general there exist several ways\footnote{For instance, for some fixed $k, l$, one can always add a term $p_k$ to $M^{(m)}_l$ and $- p_l$ to $M^{(m)}_k$ without changing $M^{(m)}$.} of writing $M^{(m)}=\sum_j p_j M^{(m)}_j$ but one looks for a set of $M^{(m)}_j$s such that condition \ref{enum:Algebra} (homogeneity of $[t^n]L_iZ$ and commutation relation \eqref{CommutationRelation}) holds which is a very much non-trivial requirement. 


\begin{proof}
The proof is based on the proof of Lemma $2.6$ in~\cite{BonzomChapuyDolega2021} which we adapt to our conditions. We proceed by induction on the coefficients of the the formal power series $Z$ in $s$. Since $[t^n]Z$ is homogeneous of degree $n$, $[t^0]Z$ is independent of the $p_i$s, hence $[t^0]L_i Z = 0$ (in fact the assumption that $[t^n]L_i F$ is homogeneous of degree $n-i$ gives that it vanishes for $n<i$). Let $n>0$ and assume that for all $n'<n$ we have $\left[t^{n'}\right]L_i Z = 0$. Since $Z$ satisfies the evolution equation, we have
\begin{equation}
    \left[ t\frac{\partial}{\partial t} - \sum_{m=1}^r t^m M^{(m)}, L_i\right] Z = \left( t\frac{\partial}{\partial t} - \sum_{m=1}^r t^m M^{(m)}\right) L_i Z.
\end{equation}
We extract the coefficient of $t^n$ and use $[t^n] t\frac{\partial F}{\partial t} = n[t^n] F$ for any formal power series $F$, as well as the induction hypothesis $[t^{n-m}]L_i Z =0$ for $m=1, \dotsc, \min(n,r)$, so that
\begin{equation} \label{eq:sfwd_exp}
[t^n]\left[ t\frac{\partial}{\partial t} - \sum_{m=1}^r t^m M^{(m)}, L_i\right] Z = n[t^n] L_iZ.
\end{equation}

Let us calculate the commutator on the LHS independently. We start with $[t\frac{\partial}{\partial t}, L_i] Z$ which is evaluated by noting that
\begin{itemize}
\item since $[t^n]Z$ is homogeneous of degree $n$, $t\frac{\partial Z}{\partial t} = \sum_{j\geq 1} jp_j \frac{\partial Z}{\partial p_j}$.
\item since $[t^n]L_i Z$ is homogeneous of degree $n-i$, $t\frac{\partial L_i Z}{\partial t} = \sum_{j\geq 1} jp_j \frac{\partial L_i Z}{\partial p_j} + iL_i Z$.
\end{itemize}
Hence,
\begin{equation}
\left[t\frac{\partial}{\partial t}, L_i\right] Z = iL_i Z + \sum_{j\geq 1} j\left[p_j\frac{\partial}{\partial p_j}, L_i\right]Z
\end{equation}
Writing $M^{(m)} = \sum_{j\geq 1} p_j M^{(m)}_j$,
\begin{equation}
\begin{aligned}
\left[ t\frac{\partial}{\partial t} - \sum_{m=1}^r t^m M^{(m)}, L_i\right] Z &= iL_i Z + \sum_{j\geq 1} \left[p_j L_j, L_i\right]Z \\
& = iL_i Z + \sum_{j\geq 1} p_j \left[ L_j, L_i\right]Z 
+ \sum_{j\geq 1} \left[ p_j, L_i\right] L_jZ
\end{aligned}
\end{equation}
The commutator $\left[ L_j, L_i\right]$ is given by the lemma assumption. As for the other commutator,
\begin{equation}
\left[ p_j, L_i\right] = -i\delta_{ij} - \sum_{m=1}^r t^m B^{(m)}_{ji}
\end{equation}
where $B^{(m)}_{ji} \coloneqq [p_j, M^{(m)}_i]$ is independent of $t$. Overall one finds
\begin{equation}
\left[ t\frac{\partial}{\partial t} - \sum_{m=1}^r t^m M^{(m)}, L_i\right] Z = t\sum_{j,k\geq 1} p_j D_{ji}^k L_k Z - \sum_{m=1}^r t^m \sum_{j\geq 1} B^{(m)}_{ji} L_jZ
\end{equation}
We now extract the coefficient of $t^n$ to get
\begin{equation}
[t^n]\left[ t\frac{\partial}{\partial t} - \sum_{m=1}^r t^m M^{(m)}, L_i\right] Z = \sum_{j,k\geq 1} p_j [t^{n-1}] D_{ji}^k L_k Z - \sum_{j\geq 1} \sum_{m=1}^r B^{(m)}_{ji} [t^{n-m}] L_jZ.
\end{equation}
Our induction hypothesis implies that $[t^{n-m}]L_k Z = 0$ for all $m= 1, \ldots, n$, and therefore the RHS is found to be zero. Equating with the RHS of \ref{eq:sfwd_exp}, we find $[t^n] L_iZ=0$.
\end{proof}


\section{$b$-deformed constellations} \label{sec:Definitions}

\subsection{$b$-weights}
A map is a graph (multiple edges are allowed) embedded in a surface such that the graph complement is a finite, disjoint union of disks which are called \emph{faces}. Given a local orientation around a vertex, there is a cyclic ordering of the edges meeting there. In a neighborhood of a vertex, a corner is the portion between two consecutive edges and it can be oriented or not.

\begin{figure}
\centering
\includegraphics[scale=0.5]{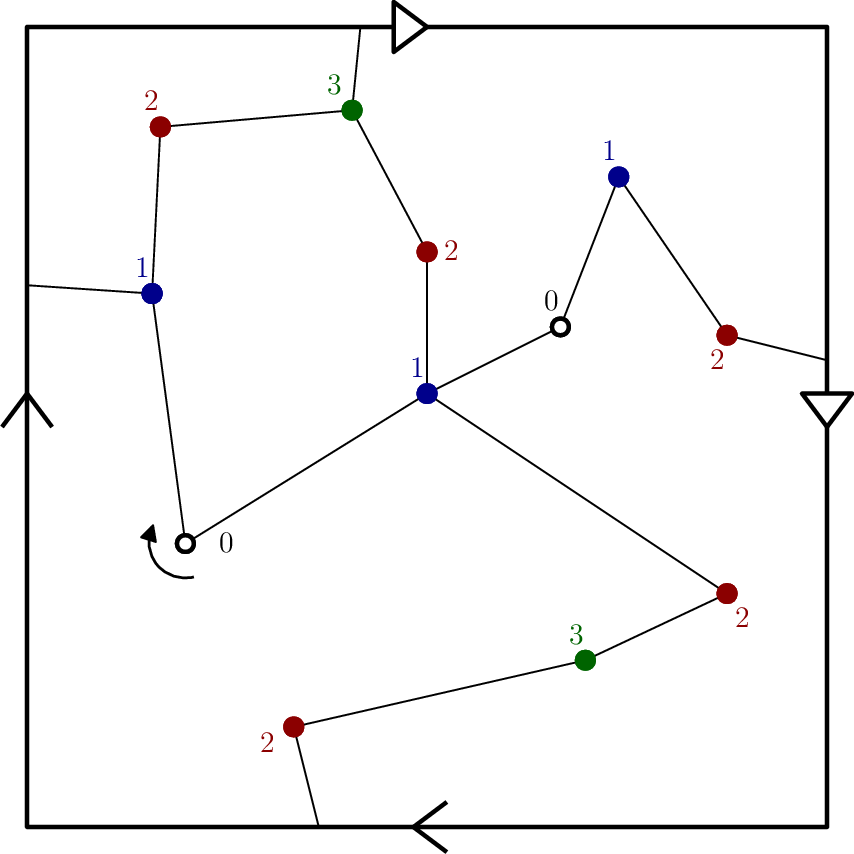}
\caption{An example of a rooted 3-constellation of size 4 on the Klein bottle.}
\label{fig:3const_Klein}
\end{figure}

Constellations are a map-like representation of ramified, orientable coverings of the sphere \cite{LandoZvonkin2004}. However, defining constellations on non-oriented surfaces is not straightforward. They have only been introduced recently in~\cite{ChapuyDolega2020} and we follow their presentation here. Let $k\geq 1$ an integer. A $k$-constellation (or simply constellation if the value of $k$ is clear or not relevant) is a map whose vertices carry a color in $\{0, \ldots, k\}$ and such that
\begin{itemize}
\item Vertices of color 0 only have neighbors of color 1,
\item Vertices of color $k$ only have neighbors of color $k-1$,
\item The cyclic order around a vertex of color $l\in\{1, \dots, k-1\}$ alternates edges connected to vertices of color $l-1$ and $l+1$.
\end{itemize}
The color of a corner is the color of the vertex it is incident to. The degree $\operatorname{deg}(f)$ of a face $f$ is its number of corners of color 0 and the size of a constellation $n(\bM)$ is its total number of corners of color 0. The degree $\operatorname{deg}(v)$ of a vertex $v$ is its number of incident edges. A \emph{root} is a distinguished, oriented corner $c$ of color 0 and a rooted constellation $(\bM,c)$ is a constellation $\bM$ equipped with a root $c$. For a (rooted or non-rooted) constellation $\bM$, further denote $F(\bM)$ the set of faces, excluding the root face if any, $V_0(\bM)$ the set of vertices of color $0$, and $v_l(\bM)$ the number of vertices of color $l\in\{1, \ldots, k\}$. An example of $3$-constellation is given in Figure~\ref{fig:3const_Klein}.

To introduce the $b$-weights, it is useful to consider the representation of maps as ribbon graphs, obtained by taking a tubular neighborhood of the embedded graph. This shows that there are sometimes two inequivalent ways of adding an edge between two corners of a map, which differ by a ``twist''. If $\bM$ is a constellation and $e$ an edge in $\bM$, we define the $b$-weight via a \emph{measure of non-orientability} $\rho_b(\bM,e)$, a function satisfying the following properties.

\begin{definition}[Measure of non-orientability in \cite{ChapuyDolega2020}] Let $\bN = \bM\setminus \{e\}$ and $c_1, c_2$ the two corners of $\bN$ onto which $e$ is attached to form $\bM$. The $b$-weight $\rho_b(\bM,e)$ must only depend on the connected component which contains $e$. Furthermore
\begin{itemize}
\item If $c_1, c_2$ are in different connected components of $\bN$, then set $\rho_b(\bM,e)=1$.
\item It they are in two different faces of the same connected component, then denote $\bM'$ the map obtained by attaching $e$ to the corners $c_1, c_2$ of $\bN$ but with a twist compared to $\bM$, and impose $\rho_b(\bM',e)+\rho_b(\bM,e) = 1+b$.
\item If $c_1$ and $c_2$ lie in the same face of $\bN$, then: either $e$ cuts the face into two in which case we set $\rho_b(\bM,e)=1$, and else we set $\rho_b(\bM,e)=b$.
\end{itemize}
\end{definition}

\begin{figure}
\centering
\includegraphics[scale=0.55]{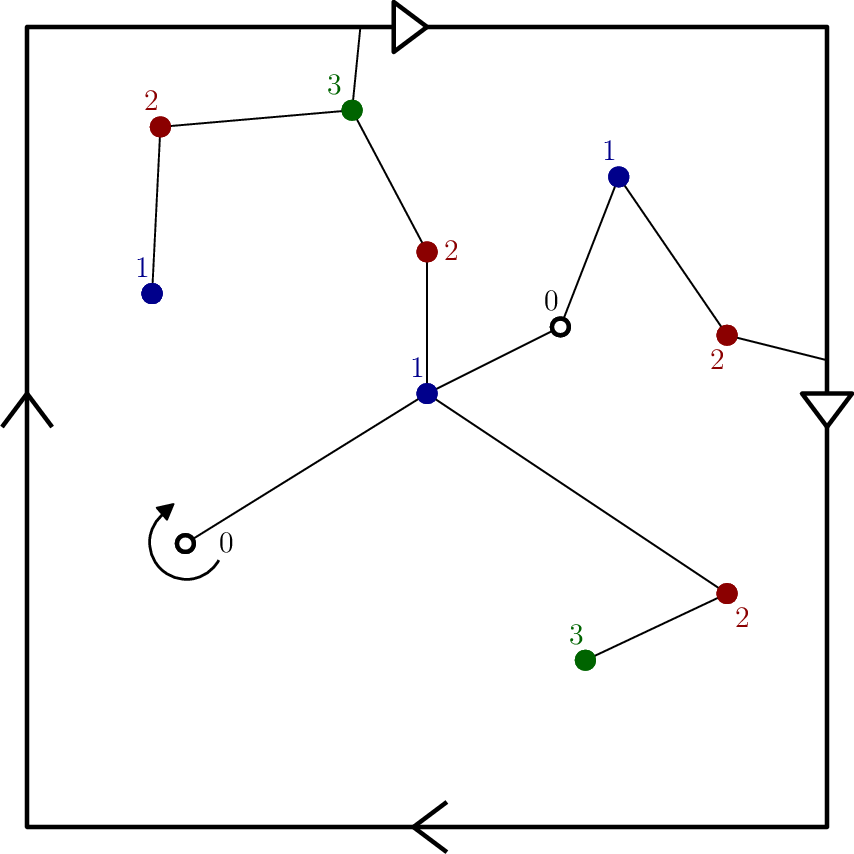}
\caption{This is the $3$-constellation obtained after deleting the right-path starting at the root in the constellation of Figure~\ref{fig:3const_Klein}. Here it is still connected and can now be embedded on the projective plane. }
\label{fig:3const_del_path}
\end{figure}

The $b$-weight of a constellation is obtained by ``collecting'' the $b$-weights of its edges in a canonical way. If $(\bM,c)$ is a rooted $k$-constellation with root corner $c$ incident to the face $f$, its \emph{right path} is the (ordered) sequence of $k$ edges $L_c=(e_1, \ldots, e_k)$ starting at $c$ and following $f$ with the orientation determined by $c$. In particular, $e_l$ connects a vertex of color $l-1$ to a vertex of color $l$, for $l=1, \ldots, k$. The $b$-weight of the right path $L_c$ is defined as
\begin{equation}
\rho_b(\bM,L_c) = \rho_b(\bM,e_1) \rho_b(\bM\setminus\{e_1\}, e_2) \dotsm \rho_b(\bM\setminus\{e_1, \ldots, e_{k-1}\}, e_k)
\end{equation}
Notice that removing from $(\bM,c)$ its right path yields another $k$-constellation (see Figure \ref{fig:3const_del_path}), which may not be connected. If so, one can still give a deterministic procedure to root each connected component~\cite{ChapuyDolega2020} but we will not need its detail here.

The algorithm of~\cite{ChapuyDolega2020} which defines the $b$-weight of $(\bM, c)$ is as follows. Say that $c$ sits at a vertex $v$ of degree $m\geq 1$, then define $(\bM_1, c_1)$ the connected component containing $v$ after removing the right path $L_c$, rooted as the corner $c_1$ inherited from $c$, and define similarly $(\bM_i, c_i)$ for $i=1, \ldots, m$ by removing the right path $L_{c_{i-1}}$, until $(\bM_m, c_m)$ which is only $v$. Set
\begin{equation}
{\rho_b}(\bM,v) = \rho_b(\bM,L_c) \rho_b(\bM_1,L_{c_1}) \dotsm \rho_b(\bM_{m-1},L_{c_m-1}).
\end{equation}

\begin{figure}
\centering
\includegraphics[scale=0.55]{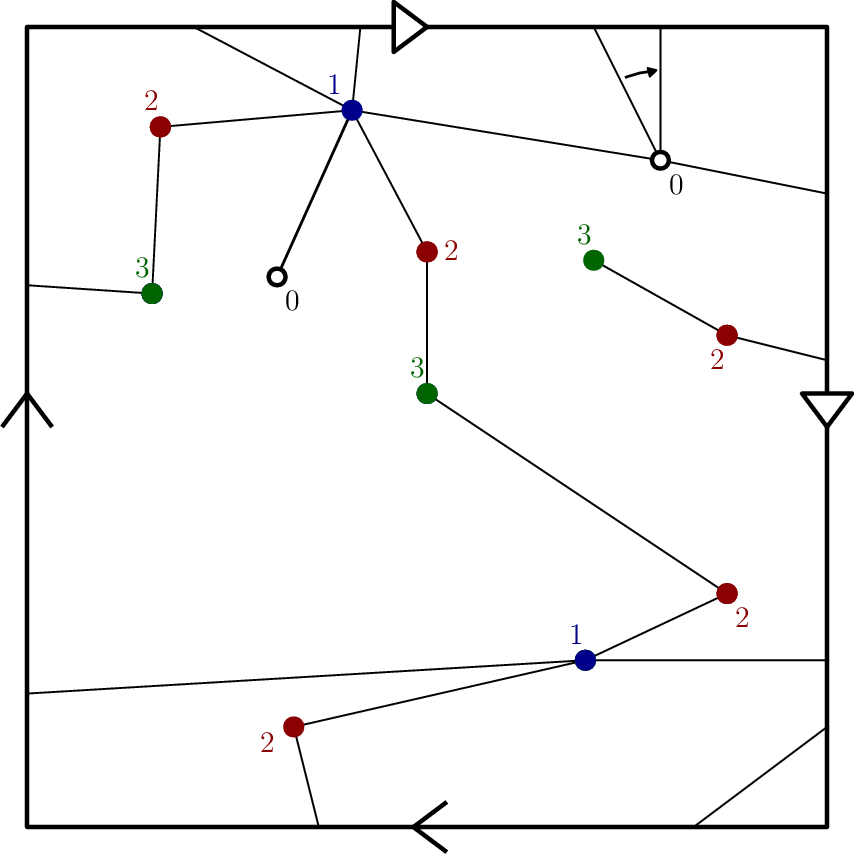}
\caption{The dual of the constellation represented in Figure~\ref{fig:3const_Klein}.}
\label{fig:3const_Klein_dual}
\end{figure}

To finish introducing the evaluation of the $b$-weights, we need a notion of dual constellations. If $\bM$ is a $k$-constellation, add a vertex of color $k+1$ for each face, connected to the corners of color $k$ that face is incident to. Then remove the vertices of color 0 and relabel the colors as $l \to k+1-l$. This gives the constellation $\tilde{\bM}$ and if $(\bM, c)$ is rooted, one can also obtain the rooted dual map $(\tilde{\bM}, \tilde{c})$, see e.g. Figure \ref{fig:3const_Klein_dual}.

Starting from $(\bM, c)$ and removing the right paths $L_c, \ldots, L_{c_{m-1}}$ as above produces some connected components $(\bM'_1,c'_1), \ldots, (\bM'_K, c'_K)$. The $b$-weight $\rho_b(\bM, c)$ is defined inductively as
\begin{equation}
\rho_b(\bM, c) = {\rho_b}(\bM,v) \rho_b(\tilde{\bM}'_1, \tilde{c}'_1) \dotsm \rho_b(\tilde{\bM}'_K, \tilde{c}'_K).
\end{equation}

\begin{figure}[!h]
\subfloat{\includegraphics[scale=0.45]{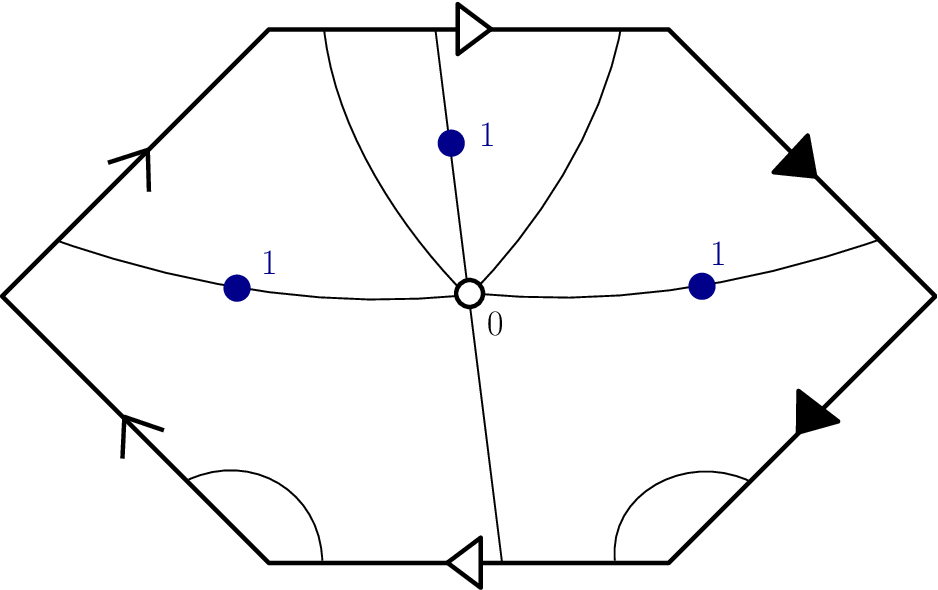}}
\hfill
\subfloat{\includegraphics[scale=0.45]{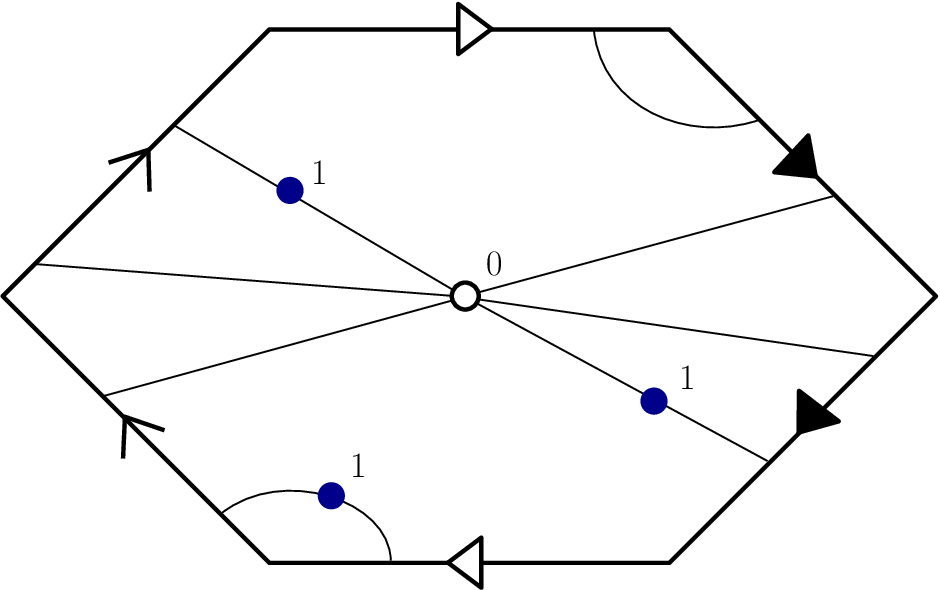}}
\caption{An example of a 2-constellation and its dual giving different $b$-weights. }
\label{fig:map_ex_bweight}
\end{figure}

Note that a constellation and its dual can have different $b$-weights. This is illustrated on Figure~\ref{fig:map_ex_bweight}: both maps have exactly one vertex of color and in both cases, there are three possible rootings, which give rise to different orders for the deletions of the right-paths. The map on the LHS has $b$-weight $b^3+b^2+b$ (each rooting gives a different $b$-weight), while the one on the RHS has $b$-weight $b^3+2b^2$ (two rootings give the same $b$-weight $b^2$).

\subsection{Generating functions}
For a rooted constellation $(\bM,c)$, denote $f_c$ the root face, and $\operatorname{deg}(f_c)$ its degree, $v_c$ the root vertex and $\operatorname{deg}(v_c)$ its degree. We consider the following sets of indeterminates 
\begin{itemize}
\item $\yy = (y_0, y_1, \ldots)$ which tracks the degree of the root face,
\item $\pp = (p_1, p_2, \ldots)$ which tracks the degrees of faces from $F(\bM)$,
\item $\qq = (q_1, q_2, \ldots)$ which tracks the degrees of vertices of color 0 from $V_0(\bM)$,
\item $u_1, \dotsc, u_k$ where $u_i$ counts the number of vertices of color $i=1, \dotsc, k$.
\end{itemize}
Then let $H^{(a, m)}_k(t,\pp,\qq, u_1, \dotsc, u_k)$ be the generating function of connected, rooted $k$-constellations with a root vertex of degree $m$ and a root face of degree $a$,
\begin{equation} 
H^{(a, m)}_k(t,\pp,\qq, u_1, \dotsc, u_k) = \sum_{\substack{(\bM,c)\\ \operatorname{deg}(f_c) = a\\ \operatorname{deg}(v_c) = m}} t^{n(\bM)} \rho_b(\bM,c) \prod_{f\in F(\bM)} p_{\operatorname{deg}(f)} \prod_{v\in V_0(\bM)} q_{\operatorname{deg}(v)} \prod_{l=1}^k u_l^{v_l(\bM)}.
\end{equation}
and $H^{(a)}_k(t,\pp,\qq, u_1, \dotsc, u_k)$ that of connected, rooted $k$-constellations with a root face of degree $a$
\begin{equation} \label{RootedMapGF}
H^{(a)}_k(t,\pp,\qq, u_1, \dotsc, u_k) = \sum_{m\geq 1} q_m t^m H^{(a, m)}(t,\pp,\qq, u_1, \dotsc, u_k).
\end{equation}
We will also need the following generating function,
\begin{equation} \label{DualGF}
\tilde{H}^{[i]}_k(t,\pp,\qq, u_1, \dotsc, u_k) = \sum_{\substack{(\bM,c)\\ \operatorname{deg}(f_c) = i}} t^{n(\bM)} \rho_b(\tilde{\bM},\tilde{c}) \prod_{f\in F(\bM)} p_{\operatorname{deg}(f)} \prod_{v\in V_0(\bM)} q_{\operatorname{deg}(v)} \prod_{l=1}^k u_l^{v_l(\bM)},
\end{equation}
which is the same as~\eqref{RootedMapGF} but calculated with the $b$-weights of the dual maps. 

\begin{proposition}\cite[Theorem 3.10]{ChapuyDolega2020}
The above generating functions satisfy
\begin{equation} \label{CombinatorialDecomposition}
H^{(a, m)}_k(t,\pp,\qq, u_1, \dotsc, u_k) = [y_a] \biggl(Y_+ \prod_{l=1}^k \Bigl(\Lambda_Y + u_l + \sum_{i,j\geq 1} y_{i+j-1} \tilde{H}_k^{[i]}(t,\pp,\qq) \frac{\partial}{\partial y_{j-1}}\Bigr)\biggr)^m y_0
\end{equation}
where
\begin{align*}
\Lambda_Y &= (1+b) \sum_{i,j\geq1} i y_{i+j-1} \frac{\partial^2}{\partial p_{i} \partial y_{j-1}} + \sum_{i,j\geq 1} p_i y_{j-1} \frac{\partial}{\partial y_{i+j-1}} + b \sum_{i\geq 1} iy_i \frac{\partial}{\partial y_i},\\
Y_+ &= \sum_{j\geq 1} y_j \frac{\partial}{\partial y_{j-1}}
\end{align*}
and if $F = \sum_{a\geq 0} y_a F_a$, then $[y_a] F \equiv F_a$ denotes the coefficient extraction with respect to $\yy$.
\end{proposition}

We explain the combinatorial origin of the terms of the equation above and refer to~\cite[Proposition 4.4 and Proposition 4.5]{ChapuyDolega2020} for details.

\begin{proof}
The operators on the RHS translate the way the $b$-weights are defined at the level of the generating function. Indeed, any constellation with a root vertex of degree $m$ can be constructed uniquely by reversing the combinatorial decomposition used to define $\rho_b(\bM,c)$ inductively. Start by creating the root vertex $v$ of color 0, equipped with an arbitrary orientation, and initially with degree 0 hence a $y_0$. For each $m'$ from 1 to $m$, we construct a new right path via the operator 
\begin{equation*}
Y_+ \prod_{l=1}^k \Bigl(\Lambda_Y + u_l + \sum_{i,j\geq 1} y_{i+j-1} \tilde{H}^{[i]}(t,\pp,\qq) \frac{\partial}{\partial y_{j-1}}\Bigr).
\end{equation*}
For each color $l=1$ to $k$, it adds an edge which changes the root face in several possible ways. 
\begin{itemize}
\item The term  $u_l$ connects the new edge to a new vertex of color $l$.
\item The term $(1+b) i y_{i+j-1} \frac{\partial^2}{\partial p_{i} \partial y_{j-1}}$ considers a root face of degree $j-1$, connects the new edge to a corner of color $l$ of the same connected component in a face of degree $i$, thus creating a root face of degree $i+j-1$. The factor $(1+b)$ is not immediate, see \cite[Proposition 4.4 and its proof]{ChapuyDolega2020}.
\item The term $p_i y_{j-1} \frac{\partial}{\partial y_{i+j-1}}$ considers a root face of degree $i+j-1$ and connects the new edge to a corner of color $l$ on that face in a way which splits it into a non-root face of degree $i$ and a root face of degree $j$.
\item The term $b iy_i \frac{\partial}{\partial y_i}$ considers a root face of degree $i$ and connects the new edge to a corner of color $l$ on that face in a way which does not split.
\item The term $y_{i+j-1} \tilde{H}^{[i]}(t,\pp,\qq) \frac{\partial}{\partial y_{j-1}}$ considers a root face of degree $j-1$ and connects the new edge to a new connected component rooted on a face of degree $i$. More precisely the new edge is attached to the corner of color $l$ on its right path, see \cite[Proposition 4.5 and its proof]{ChapuyDolega2020}.
\end{itemize}
After adding a right path, we update the fact that we have a new $k$-constellation whose root face has one more corner of color 0, hence the action of $Y_+$.
\end{proof}

\subsection{Deriving the constraints by applying Lemma~\ref{lemma:constr_from_ev}} 
\subsubsection{The evolution equation} A constellation is \emph{labeled} if its $n$ corners of color 0 are labeled and oriented. We then form the generating function of connected labeled $k$-constellations
\begin{equation}
H_k(t,\pp,\qq, u_1, \dotsc, u_k) = \sum_{\text{labeled $\bM$}} \frac{t^{n(\bM)} \rho_b(\bM)}{2^{n(\bM)-1} n(\bM)!} \prod_{f\in F(\bM)} p_{\operatorname{deg}(f)} \prod_{v\in V_0(\bM)} q_{\operatorname{deg}(v)} \prod_{l=1}^k u_l^{v_l(\bM)},
\end{equation}
where $\rho_b(\bM)$ is defined by \eqref{AverageWeight}. Finally, let 
\begin{equation}
\tau_k(t,\pp,\qq, u_1, \dotsc, u_k) = \exp \frac{H_k(t,\pp,\qq, u_1, \dotsc, u_k)}{1+b}
\end{equation}
be the generating of non-necessarily connected labeled $k$-constellations. It was proved in~\cite{ChapuyDolega2020} using heavy algebraic manipulations that it satisfies the following evolution equation.
\begin{theorem}\cite{ChapuyDolega2020}
The generating function $\tau_k(t,\pp,\qq)$ satisfies
\begin{equation}
\label{eq:decomp_eq_const}
\frac{\partial \tau_k(t,\pp,\qq, u_1, \dotsc, u_k)}{\partial t} = \sum_{m\geq 1} q_m t^{m-1} \Theta_Y \Bigl(Y_+ \prod_{l=1}^k \bigl(\Lambda_Y+u_l\bigr)\Bigr)^m \frac{y_0}{1+b} \tau_k(t,\pp,\qq, u_1, \dotsc, u_k),
\end{equation}
with $\Theta_Y = \sum\limits_{i\geq 1} p_i [y_i]$.
\end{theorem}
Notice that $\Bigl(Y_+ \prod_{l=1}^k \bigl(\Lambda_Y+u_l\bigr)\Bigr)^m \frac{y_0}{1+b}$ is linear in $\yy$ so the action of $\Theta_Y$ is to replace the variables marking the root face with the ordinary variables for faces.

We wish to apply Lemma \ref{lemma:constr_from_ev} to the series $\tau_k$. Four conditions have to be checked. Condition \ref{enum:Homogeneity} is clear: $[t^n]\tau_k$ is a homogeneous polynomial of degree $n$ in the $p_i$s (with the rule that $p_i$ has degree $i$). The above evolution equation \eqref{eq:decomp_eq_const} gives the evolution operators
\begin{equation} \label{Mkm}
M^{(k,m)} \coloneqq q_m\Theta_Y \Bigl(Y_+ \prod_{l=1}^k \bigl(\Lambda_Y+u_l\bigr)\Bigr)^m \frac{y_0}{1+b},
\end{equation}
for $m\geq 1$. It is differential in the $p_i$s, linear in $q_m$ and symmetric in $u_1, \ldots, u_k$. 

Obviously one can insert a bound $r$ such that $q_{m}=0$ for $m>r$. Formally, let
\begin{equation}
\tau_{k,r} \equiv \tau_{k,r}(t, \pp, q_1, \dotsc, q_r, u_1, \dotsc, u_k) \coloneqq \tau_k(t, \pp, \qq, u_1, \dots, u_k)_{|q_m = 0\ \forall m>r},
\end{equation}
then \eqref{eq:decomp_eq_const} reads $t\frac{\partial\tau_{k,r}}{\partial t} = \sum_{m=1}^r t^m M^{(k,m)} \tau_{k,r}$, where $M^{(k,1)}, \ldots, M^{(k,r)}$ are independent of $t$. This is condition \ref{enum:Evolution}. 

We are left with condition \ref{enum:Algebra} which we are not able to check beyond the models mentioned in the Introduction. Let us introduce the candidate operators $M^{(k,m)}_i$. For $m,k,i\geq1$, let
\begin{equation}
    \label{eq:def_modes}
     M^{(k,m)}_i \coloneqq \left[y_i\right] \left(Y_+\prod\limits_{c=1}^{k} (\Lambda_Y+u_c)\right)^m \frac{y_0}{1+b}
\end{equation}
Then the evolution operator on the RHS of \eqref{eq:decomp_eq_const} reads
\begin{equation}
\sum_{m=1}^r q_m t^{m-1} \Theta_Y \Bigl(Y_+ \prod_{l=1}^k \bigl(\Lambda_Y+u_l\bigr)\Bigr)^m \frac{y_0}{1+b} = \sum_{i\geq 1} \sum_{m=1}^r t^{m-1} q_m p_i M^{(k,m)}_i.
\end{equation}
We want to check that the evolution operators satisfy condition~\ref{enum:Algebra} by computing their commutators. The complexity of this computation essentially depends on the highest power of the operator $\Lambda_Y$ which appears in the expression of $M^{(k,m)}_i$. In the model with evolution operator $M^{(k,m)}_i$, the highest power of $\Lambda_Y$ is $km$. Accordingly, we say that the constellation model where $q_j = \delta_{j,m}$ for $j\geq 1$, with evolution operator $M^{(k,m)}$, has order $km$. 

Using the fact that $t\frac{\partial\tau_{k,r}}{\partial t} = \sum_{i\geq 1} ip_i \frac{\partial \tau_{k,r}}{\partial p_i}$, the evolution equation can then be rewritten as
\begin{equation}
\sum_{i\geq 1} p_i L^{(k,r)}_i \tau_{k,r} = 0, \quad\text{with} \quad
L^{(k,r)}_i \coloneqq -i\frac{\partial}{\partial p_i} + \sum_{m=1}^r t^m q_m M^{(k,m)}_i.
\end{equation}

We now specialize to cases where we have been able to prove condition \ref{enum:Algebra} and thus extract constraints. In the remaining we use the notation $p_i^* \equiv i\frac{\partial}{\partial p_i}$ for $i\geq 1$ and the convention $p_j^*\equiv 0$ for $j\leq 0$.

\subsubsection{Bipartite maps} They correspond $k=2$ and $q_m = \delta_{m,1}$ (so that $r=1$). This way, the vertices of color 0 are leaves and can be erased, while the vertices of colors 1 and 2 are counted by the variables $u_1, u_2$. Hence
\begin{equation*}
\begin{aligned}
L_i^{\text{bip}} &\coloneqq L^{(2,1)}_{i|q_m = \delta_{m,1}} = -p_i^* + [y_i] t Y_+ \bigl(\Lambda_Y+u_1\bigr)\bigl(\Lambda_Y+u_2\bigr) \frac{y_0}{1+b}\\
&= -p_i^* + (1+b)t \sum_{l,m\geq 1\atop l+m=i-1} p_l^* p_m^* + t\sum_{l\geq 1} p_l p_{l+i-1}^* + t\bigl(b (i-1) + u_1 +u_2\bigr) p_{i-1}^* + t\frac{u_1u_2}{1+b} \delta_{i,1}
\end{aligned}
\end{equation*}

A direct computation shows that these operators form a (shifted) Virasoro algebra.
\begin{theorem} \label{thm:BipMaps}
They satisfy for all $i\geq 1$
\begin{equation}
[L_i^{\text{bip}}, L_j^{\text{bip}}] = t(i-j)L_{i+j-1}^{\text{bip}}.
\end{equation}
\end{theorem}

While this result is technically new, since the notion of $b$-deformed bipartite maps was only introduced in \cite{ChapuyDolega2020}, it should not come as a surprise. In the case of general maps (not necessarily bipartite; a case we include below as a special case of Theorem \ref{thm:2ndModel}), \cite{Okounkov1997} provides a relation between the series of Goulden and Jackson \eqref{GouldenJacksonSeries} and the $\beta$-ensemble of the 1-matrix integral. Then \cite{AdlervanMoerbeke2001} derives general Virasoro constraints for this type of integrals. Notice that the interpretation using $b$-weights only comes later with \cite{LaCroix2009}. This is for general maps, but there is no surprise that it works the same for bipartite maps since they both satisfy the (half-)Virasoro algebra in the orientable case\footnote{Going from general maps to bipartite maps translates in the world of matrix integrals as going from the Gaussian ensemble to the Laguerre ensemble, which both satisfy (half-)Virasoro constraints.}.

The method of~\cite{AdlervanMoerbeke2001} and the one presented here are very different. In~\cite{AdlervanMoerbeke2001}, the constraints are identified via the connection with the $\beta$-ensemble of matrix integrals (without combinatorial interpretation as maps), while we rely on the combinatorial techniques of~\cite{ChapuyDolega2020}.

\subsubsection{3-constellations} They correspond to $k=3$ and $q_m = \delta_{m,1}$ (again the vertices of color 0 are leaves which can be erased while the vertices of colors 1, 2, 3 are counted by the variables $u_1, u_2, u_3$) hence
\begin{equation*}
L_i^{\text{3-const.}} \coloneqq L^{(3,1)}_{i|q_m = \delta_{m,1}} = -p_i^* + [y_i] t Y_+ \bigl(\Lambda_Y+u_1\bigr)\bigl(\Lambda_Y+u_2\bigr)\bigl(\Lambda_Y+u_3\bigr) \frac{y_0}{1+b}.
\end{equation*}
Explicitly,
\begin{multline*}
L_i^{\text{3-const.}} = -p_i^* + (1+b)^2t \sum_{l,m,n\geq 1\atop l+m+n=i-1} p_l^* p_m^* p_n^* + (1+b)t \Bigl(\sum_{m,n\geq 1\atop m+n\geq i} + \sum_{m=1}^{i-1} \sum_{n\geq 1\atop m+n\geq i}\Bigr) p_{n+m-i+1} p_m^* p_n^*\\
+ t\sum_{m,n\geq 1} p_n p_m p_{n+m+i-1}^* + (1+b) t\frac{i(i-1)}{2} p_{i-1}^* + b(1+b) t\frac{3}{2}(i-1)\sum_{m,n\geq 1\atop m+n=i-1} p_n^* p_m^*\\
+ b^2t (i-1)^2 p_{i-1}^* + b \sum_{n\geq 1} (n+2i-2) p_n p_{n+i-1}^*\\
+ t (u_1 + u_2 + u_3) \Bigl( (1+b) \sum_{m,n\geq 1\atop m+n=i-1} p_n^* p_m^* + \sum_{n\geq 1} p_n p_{n+i-1}^* + b(i-1)p_{i-1}^*\Bigr)\\
+ t (u_1 u_2 + u_1 u_3 + u_2 u_3) p_{i-1}^* + t\frac{u_1 u_2 u_3}{1+b} \delta_{i,1}.
\end{multline*}

\begin{theorem} \label{thm:3-constellations}
They satisfy for all $i\geq 1$
\begin{multline}
[L_i^{\text{3-const.}}, L_j^{\text{3-const.}}] = 2t(i-j) \sum_{n\geq 1} p_n L^{\text{3-const.}}_{i+j+n-1} + t b(i-j)(i+j-2) L^{\text{3-const.}}_{i+j-1}\\
+ (1+b)t\Bigl(3(i-j) \sum_{n=1}^{\mu-1} + \operatorname{sgn}(i-j) \sum_{n=\mu}^{M-1} (2M-2n-\mu-1)\Bigr) p_n^* L_{i+j-1-n}^{\text{3-const.}} \\
+ t(u_1 + u_2 + u_3)(i-j) L^{\text{3-const.}}_{i+j-1}
\end{multline}
with $\mu\coloneqq \min(i,j)$ and $M\coloneqq\max(i,j)$. 
\end{theorem}
Notice that there are two places with an explicit $b$-dependence. The coefficient $(1+b)$ can be re-absorbed by defining $p_i^\perp = (1+b)p_i^*$ which is in fact natural from the point of view of the theory of symmetric functions\footnote{Then $p_i^\perp$ is the dual function to $p_i$ for the $\alpha$-deformed scalar product, itself defined by $\langle p_\lambda, p_\nu\rangle_\alpha = \alpha^{l(\lambda)}\delta_{\lambda, \nu} z_\lambda$ with $\alpha=1+b$ for us.}. However it is not possible to get rid of the other appearance of $b$ in the expression of commutator above.

\subsubsection{Bipartite maps with black vertices of degree less than or equal to 3} They correspond to $k=1$ (so that vertices of color 1 are counted by a variable $u$ and there are no colors $c>1$) and $r=3$, i.e. $q_m = 0$ for $m>3$, so that $q_1, q_2, q_3$ control the degrees of the vertices of color 0. We have
\begin{equation}
L^{\text{bip}\leq 3}_i \coloneqq L_i^{(1,3)} = -p_i^* + \sum_{m=1}^3 q_m t^m M^{(1,m)}_i
\end{equation}
Explicitly,
\begin{equation} \label{M1m}
\begin{aligned}
M^{(1,1)}_i &= p_{i-1}^* + u\delta_{i,1},\\ 
M^{(1,2)}_i &= (1+b)t^2\sum_{n,m\geq1\atop n+m=i-2} p_n^* p_m^* + t^2\sum_{n\geq 1} p_n p_{n+i-2}^* + t^2 (b(i-1)+2u) p_{i-2}^*\\ &+ \frac{t^2 u}{1+b}p_1 \delta_{i,1} + \frac{t^2u(b+u)}{1+b}\delta_{i,2}
\end{aligned}
\end{equation}
while the explicit expression of $M^{(1,3)}_i$ in terms of $p_i$s and $p_i^*$s is quite cumbersome and we do not write it. 

\begin{theorem} \label{thm:2ndModel}
They satisfy for all $i,j \geq 1$
\begin{multline}
[L_i^{\text{bip}\leq 3}, L_j^{\text{bip}\leq 3}] = 2t^3q_3(i-j) \sum_{n\geq 1} p_n L^{\text{bip}\leq 3}_{i+j+n-3} + t^3 q_3b(i-j)(i+j-3) L^{\text{bip}\leq 3}_{i+j-3}\\
+ (1+b)t^3 q_3\Bigl(3(i-j) \sum_{n=1}^{\mu-2} + \operatorname{sgn}(i-j) \sum_{n=\mu-1}^{M-2} (2M-2n-\mu-3)\Bigr) p_n^* L_{i+j-3-n}^{\text{bip}\leq 3} \\
+ 3t^3q_3u(i-j) L^{\text{bip}\leq 3}_{i+j-3} + t^2 q_2 (i-j)L^{\text{bip}\leq 3}_{i+j-2}.
\end{multline}
\end{theorem}
This algebra already appeared at $b=0$ in \cite{MarshakovMironovMorozov}. This model also contains $b$-deformed general maps, obtained by setting $q_1=q_3=0$. In this case, Theorem~\ref{thm:2ndModel} gives that $b$-deformed maps also have Virasoro constraints. As discussed above, this was known from \cite{AdlervanMoerbeke2001} in the context of matrix integrals, and using \cite{Okounkov1997} and \cite{LaCroix2009} to make the connection with maps. 

\subsection{Corollaries \ref{thm:Rooting} and \ref{thm:Duality}} Both corollaries follow from interpreting the constraints $L_i^{(k,r)}\tau_{k,r} = 0$ combinatorially on \emph{connected} constellations. In the models we consider we found
\begin{equation}
i\frac{\partial \tau_{k,r}}{\partial p_i} = \sum_{m=1}^r t^{m-1} q_m [y_i] \left(Y_+ \prod_{l=1}^k (\Lambda_Y+ u_l)\right)^m \frac{y_0}{1+b}\tau_{k,r}.
\end{equation}
We will insert $\tau_{k,r} = e^{H_{k,r}/(1+b)}$ to try and prove Corollary \ref{thm:Rooting}, i.e. $i\frac{\partial H_{k,r}}{\partial p_i} = H^{(i)}_{k,r}$. We use
\begin{equation}
\Lambda_Y \tau_{k,r} = \tau_{k,r} \left(\Lambda_Y + \sum_{i,j\geq 1} y_{i+j-1} \frac{i\partial H_{k,r}}{\partial p_i} \frac{\partial}{\partial y_{j-1}}\right)
\end{equation}
and thus
\begin{multline} \label{ConstraintRewritten}
i\frac{\partial H_{k,r}}{\partial p_i} 
= \sum_{m=1}^r t^{m-1} q_m [y_i] \\ \left(Y_+ \prod_{l=1}^k \left(\Lambda_Y+ u_l + \sum_{i,j\geq 1} y_{i+j-1} \frac{i\partial H_{k,r}}{\partial p_i} \frac{\partial}{\partial y_{j-1}}\right)\right)^m y_0.
\end{multline}

From \cite[Corollary 5.9]{ChapuyDolega2020} we know that
\begin{equation} \label{VertexRooting}
m\frac{\partial H_k(t,\pp,\qq, u_1, \dotsc, u_k)}{\partial q_m} = H^{[m]}_k(t,\pp,\qq, u_1, \dotsc, u_k)
\end{equation}
where $H^{[m]}_k(t,\pp,\qq, u_1, \dotsc, u_k)$ is the generating function of connected $k$-constellations with a root vertex of degree $m$ (and no restrictions on $\qq$ for this equation to make sense). Let $\pi$ be the operator which exchanges $\pp$ and $\qq$. By applying duality, it comes that $\pi H^{[m]}_k = \tilde{H}^{[m]}_k$ where the latter is defined in \eqref{DualGF}. It is the generating function of connected $k$-constellations with a root face of degree $m$ and whose $b$-weights are those of their dual maps. Therefore, applying $\pi$ to \eqref{VertexRooting} gives
\begin{equation} \label{DualRooting}
i\frac{\partial H_k(t,\pp,\qq, u_1, \dotsc, u_k)}{\partial p_i} = \tilde{H}^{[i]}_k(t,\pp,\qq, u_1, \dotsc, u_k).
\end{equation}
The RHS of \eqref{ConstraintRewritten} becomes exactly that of \eqref{CombinatorialDecomposition} after evaluation on $q_m = 0$ for $m>r$, i.e. the combinatorial decomposition of $H^{(i)}_{k,r}$. We thus get Corollary \ref{thm:Rooting}. Corollary \ref{thm:Duality} is a simple rewriting, by replacing the LHS using \eqref{DualRooting} (again after evaluation on $q_m = 0$ for $m>r$).

\section{Computation of the commutation relations for bipartite maps and 3-constellations} \label{sec:3Constellations}

Let $\alpha\in\mathbb{R}$ and recall the notation $p_i^* = i\frac{\partial}{\partial p_i}$. Let the ``$b$-deformed currents'' $J^{(b)}_i$ for $i\in\mathbb{Z}$ be
\begin{equation}
    \label{eq:currents}
    J^{(b)}_i = \begin{cases} p_{-i} \qquad &\text{for $i<0$,} \\ (1+b)p_i^* \qquad &\text{for $i>0$,}\\ \alpha \qquad  & \text{for $i=0$,}\end{cases}
\end{equation}
which satisfy the commutation relation
\begin{equation}
    \left[J^{(b)}_i,J^{(b)}_j\right] = (1+b)i\delta_{i,-j},
\end{equation}
It is a representation of the Heisenberg algebra where $\alpha$ is a central element called the charge. {\bf In the rest of this section we set $\alpha=0$.}

\subsection{Proof of Theorems \ref{thm:BipMaps}, \ref{thm:3-constellations}}
Recall that the operators 
\begin{equation*}
M^{(k,m)}_i = [y_i] \left(Y_+ \prod_{c=1}^k (\Lambda_Y+u_c)\right)^m \frac{y_0}{1+b}
\end{equation*}
are symmetric polynomials in $u_1, \ldots, u_k$. We will use their expansion onto the elementary symmetric polynomials $e_{p}(u_1,\ldots,u_{k})$ for $0 \leq p \leq k-1$ to organize our calculation. For $i\geq 1$ and integers $s_1, \dotsc, s_m\geq 0$, let
\begin{equation}
A^m_i(s_1, \ldots, s_m) \coloneqq [y_i] \prod_{l=1}^m Y_+\Lambda_Y^{s_l}\ \frac{y_0}{1+b},
\end{equation}
so that
\begin{equation}
M^{(k,m)}_i = \sum_{s_1, \ldots, s_m =0}^k A^m_i(s_1, \ldots, s_m) \prod_{l=1}^m e_{k-s_l}(u_1, \ldots, u_k).
\end{equation}

We will focus on the following case, relevant for $k$-constellations with $q_m = \delta_{m,1}$,
\begin{equation}
A_i(s)\equiv A_i^1(s) = [y_i]Y_+\Lambda_Y^s \frac{y_0}{1+b},
\end{equation}
for $i\geq 1, s\geq 0$. 
We form the candidate constraints
\begin{equation}
L^{(k)}_i \equiv L^{(k,1)}_{i|q_m = \delta_{m,1}} = -p_i^* + t M^{(k,1)}_i = -p_i^* + t \sum_{s=0}^k e_{k-s}(u_1, \dotsc, u_k) A_i(s),
\end{equation}
for which the following proposition provides a recursion.
\begin{proposition}
One has $A_i(0) = \frac{\delta_{i,1}}{1+b}$ for all $i\geq 1$. Then for all $i\geq1, s\geq 0$
\begin{equation} \label{eq:mod_rec_no_ui}
    A_i(s+1) = \sum\limits_{\substack{n \geq 1 }} J^{(b)}_{i-n} A_n(s) + b(i-1) A_i(s) 
\end{equation}
\end{proposition}

\begin{proof}
Let $Y_i$, $i\in\mathbb{Z}$, be the operator defined as $Y_i y_j = y_{i+j} \delta_{i+j\geq 0}$. Then,
\begin{equation} \label{eq:def_lambda_J}
\Lambda_Y = \sum_{i\in\mathbb{Z}} Y_i J^{(b)}_i + b \sum\limits_{i \geq 0} iy_i [y_i].
\end{equation}
By definition, $A_i(0) = [y_i] Y_+\frac{y_0}{1+b} = \frac{\delta_{i,1}}{1+b}$ for all $i\geq 0$. 
Then
\begin{equation}
    A_i(s+1) = \left[y_i\right] Y_+ \Lambda_Y^{s+1} \frac{y_0}{1+b} = \left[y_{i-1}\right] \Lambda_Y \sum\limits_{j \geq 1} y_{j-1} \underbrace{\left[y_{j-1}\right]  \Lambda_Y^{s} \frac{y_0}{1+b}}_{A_j(s)}
\end{equation}
and using~\eqref{eq:def_lambda_J} gives~\eqref{eq:mod_rec_no_ui}.
\end{proof}

In the rest of this section, we compute the commutators for bipartite maps and $3$-constellations, i.e. for $k\leq 3$. 
Our strategy will be based on the following lemmas which reconstruct the commutators of the constraints from those of the operators $A_i(s)$.

\begin{lemma} \label{lemma:SimplifiedCommutators}
Assume that there exist differential operators $(D_{ij, l}(s))_{i,j,l\geq 1, s\geq 0}$ in the variables $p_i$s such that
\begin{equation} \label{Commutatorp*A}
[p_i^*, A_j(s)] - [p_j^*, A_i(s)] = \sum_{l\geq 1} D_{ij, l}(s) p_l^*
\end{equation}
and
\begin{equation} \label{eq:SimplifiedCommutators}
\begin{aligned}
    \left[A_i(s), A_j(s)\right] &= \sum_{l\geq 1} D_{ij,l}(s) A_l(s), \\
    \left[A_i(s), A_j(s')\right] -  \left[A_j(s), A_i(s')\right] &= \sum_{l\geq 1} D_{ij,l}(s') A_l(s) + D_{ij,l}(s) A_l(s').
\end{aligned}
\end{equation}
Then the operators $L^{(k)}_i = -p_i^* + t \sum_{s=0}^k e_{k-s}(u_1, \dotsc, u_k) A_i(s)$ satisfy
\begin{equation}
\label{eq:CommutConstr}
\left[L^{(k)}_i, L^{(k)}_j\right] = t \sum_{l\geq 1} D_{ij, l}^{(k)} L_l^{(k)},
\end{equation}
with
\begin{equation}
D_{ij, l}^{(k)} = \sum_{s=0}^k e_{k-s}(u_1,\ldots,u_k) D_{ij, l}(s).
\end{equation}
\end{lemma}

\begin{proof}
Recall that $M^{(k,1)}_i = \sum_{s=0}^k e_{k-s}(u_1,\ldots,u_k) A_i(s)$, which combined with~\eqref{eq:SimplifiedCommutators}, gives directly
\begin{equation}
\left[M^{(k,1)}_i, M^{(k,1)}_j\right] = \sum_{l\geq 1} D_{ij, l}^{(k)} M_l^{(k)}.
\end{equation}
In addition, the constraints read $L^{(k,1)}_{i|q_m = \delta_{m,1}} = -p_i^* + tM^{(k,1)}_i$, hence
\begin{equation}
[L^{(k,1)}_{i|q_m = \delta_{m,1}}, L^{(k,1)}_{i|q_m = \delta_{m,1}}] = -t\left( \left[p_i^*,M^{(k,1)}_j \right] + \left[M^{(k,1)}_i, p_j^* \right] \right) 
+ t^2\left[M^{(k,1)}_i, M^{(k,1)}_j\right].
\end{equation}
We use \eqref{Commutatorp*A} and get directly~\eqref{eq:CommutConstr}.
\end{proof}

As it turns out, for any set of operators $A_i(s)$, the equations \eqref{eq:SimplifiedCommutators} always admit the solution $D_{ij, l}(s) = \delta_{j,l} A_i(s) - \delta_{i,l} A_j(s)$, and therefore any operators of the form $M^{(k,1)}_i = \sum_{s=0}^k e_{k-s}(u_1,\ldots,u_k) A_i(s)$ satisfy $[M^{(k,1)}_i, M^{(k,1)}_j] = \sum_{l\geq 1} D_{ij, l}^{(k)} M_l^{(k)}$. The non-triviality of Lemma \ref{lemma:SimplifiedCommutators} comes the fact that the constraints read $L^{(k)}_i = -p_i^* + tM^{(k,1)}_i$, and that we assume that the coefficients $D_{ij, l}(s)$ satisfy Equation \eqref{Commutatorp*A}. For the same reason, the assumption that $[L_i, L_j] = t \sum_{k\geq 1} D_{ij,k} L_k$ in the Lemma \ref{lemma:constr_from_ev} would be trivial without the factor $t$ on the RHS, since $D_{ij,k}$ could be a combination of the $L_i$s themselves.

Our strategy is to prove the relations \eqref{Commutatorp*A} and \eqref{eq:SimplifiedCommutators}, from which the commutators~\eqref{eq:CommutConstr} and condition~\ref{enum:Algebra} follow.

\begin{proposition} \label{thm:SimplifiedCommutators}
Let $\mu = \min(i,j)$ and $M = \max(i,j)$. For $0\leq s, s'\leq 3$,
\begin{subequations} \label{SimplifiedCommutators}
\begin{alignat}{1}
    \left[A_i(s), A_j(s)\right] &= \sum_{l\geq 1} D_{ij, l}(s) A_l(s),\label{CommutatorA} \\
    \left[A_i(s), A_j(s')\right] -  \left[A_j(s), A_i(s')\right] &= \sum_{l\geq 1} D_{ij,l}(s') A_l(s) + D_{ij,l}(s) A_l(s'), \label{CommutatorAA}\\
    \left[p_i^*, A_j(s)\right] -  \left[p_j^*, A_i(s)\right] &= \sum_{l\geq 1} D_{ij,l}(s) p_l^* \label{CommutatorpA}
\end{alignat}
\end{subequations}
with $D_{ij, l}(0) = D_{ij, l}(1) = 0$ and $D_{ij, l}(2) = (i-j) \delta_{l,i+j-1}$ and
\begin{multline}
D_{ij, l}(3) = (i-j) (2\delta_{l\geq M} + \delta_{M\leq l\leq i+j-1}) J^{(b)}_{i+j-1-l} \\
+ \operatorname{sgn}(i-j)(2l-3\mu+1) \delta_{\mu\leq l\leq M-1} J^{(b)}_{i+j-1-l}
+ b(i-j)(i+j-2)\delta_{l, i+j-1}.
\end{multline}
\end{proposition}

This proposition shows that the assumptions of the Lemma \ref{lemma:SimplifiedCommutators} are satisfied by $k$-constellations for $k\leq 3$, and thus immediately implies Theorem~\ref{thm:BipMaps} for bipartite maps and Theorem~\ref{thm:3-constellations} for 3-constellations.
\begin{proof}[Proof of Theorem \ref{thm:BipMaps} for bipartite maps and Theorem \ref{thm:3-constellations} for 3-constellations]
One expands the constraints as
\begin{equation*}
\begin{aligned}
L^{\text{bip}}_i &= -p_i^* + t(A_i(2) + e_1(u_1, u_2) A_i(1) + e_2(u_1, u_2))\\
L^{\text{3-const.}}_i &= -p_i^* + t(A_i(3) + e_1(u_1, u_2, u_3) A_i(2) + e_2(u_1, u_2, u_3) A_i(1) + e_3(u_1, u_2, u_3))
\end{aligned}
\end{equation*}
It is then straightforward to feed Proposition \ref{thm:SimplifiedCommutators} to Lemma \ref{lemma:SimplifiedCommutators} and compute for bipartite maps
\begin{equation*}
[L^{\text{bip}}_i, L^{\text{bip}}_j] = t\sum_{l\geq 1} D_{ij, \ell}^{(2)} L^{\text{bip}}_\ell
\end{equation*}
with $D_{ij, \ell}^{(2)} = D_{ij,\ell}(2) + e_1(u_1, u_2) D_{ij, \ell}(1) + e_2(u_1, u_2) D_{ij, \ell}(0) = (i-j)$. Similarly for 3-constellations,
\begin{equation*}
[L^{\text{3-const.}}_i, L^{\text{3-const.}}_j] = t\sum_{\ell\geq 1} D_{ij, \ell}^{(3)} L^{\text{3-const.}}_\ell
\end{equation*}
with $D_{ij, \ell}^{(3)} = D_{ij,\ell}(3) + e_1(u_1, u_2, u_3) D_{ij, \ell}(2) + e_2(u_1, u_2, u_3) D_{ij, \ell}(1) + e_3(u_1, u_2, u_3) D_{ij, \ell}(0)$.
\end{proof}


%
%


This section presents the detailed computation, and thereby the proof, of Proposition~\ref{thm:SimplifiedCommutators}. It is sufficient to check the equations \eqref{CommutatorA}, \eqref{CommutatorAA}, \eqref{CommutatorpA} for all $s, s'\leq 3$.

We will often encounter quantities depending on two integers $i, j$, which are variations of $\sum_{n\geq i} (i-n) O_{i+j, n}$, up to index shifts, where $O_{i+j, n}$ is some operator. These quantities are then antisymmetrized. In those circumstances, we use
\begin{equation}
\sum_{n\geq i} (i-n) O_{i+j, n} - \sum_{n\geq j} (j-n) O_{i+j, n} = (i-j) \sum_{n\geq M} O_{i+j, n} - \sgn(i-j) \sum_{n=\mu}^{M-1} (\mu-n) O_{i+j, n}.
\end{equation}

In terms of notation, we use the Kronecker delta $\delta_{i,j} = 1$ if $i=j$ and 0 else, and we also use the notation $\delta_{i\leq j}$ to denote the step function which is 1 if $i\leq j$ and 0 else.

\subsection{Computation of $\left[ A_i(s), A_j(s') \right]$ for $s, s' =0, 1$}
We have $A_i(0) = \frac{\delta_{i,1}}{1+b}$ and $A_i(1) = \frac{J^{(b)}_{i-1}}{1+b}$ for $i\geq 1$. Therefore
\begin{equation*}
[A_i(0), A_j(0)] = [A_i(0), A_j(1)] = [A_i(1), A_j(1)] = [p_i^*, A_j(0)] = [p_i^*, A_j(1)] = 0
\end{equation*}
which are all compatible with $D_{ij, l}(0) = D_{ij, l}(1) = 0$ indeed. For later purposes, we also give for $j\in\mathbb{Z}$, $[J_j^{(b)}, A_i(1)] = j \delta_{i+j,1}$.

\subsection{Computation of $\left[A_i(s), A_j(s') \right]$ for $s, s' \leq 2$} \label{sec:A2A2}
We use the relation~\eqref{eq:mod_rec_no_ui} to compute the various commutators involving $A_i(2)$, starting with
\begin{equation}
\begin{aligned}
    \left[J_j^{(b)} , A_i(2) \right] &= \sum_{n\geq 1} J_{i-n}^{(b)} [J^{(b)}_j, A_n(1)] + [J^{(b)}_j, J^{(b)}_{i-n}]A_n(1) + b(i-1)[J^{(b)}_j, A_i(1)]\\
    &= jJ^{(b)}_{i+j-1}\left( \delta_{j \leq 0} + \delta_{i+j \geq 1} \right) - b(i-1)^2 \delta_{i+j,1}.
\end{aligned}
\end{equation}
A direct application gives $\left[A_i(1) , A_j(2) \right] = (i-1) A_{i+j-1}(1)$, and $\left[p_i^* , A_j(2) \right] = i p_{i+j-1}^*$, then
\begin{equation}
\begin{aligned}
\left[A_i(1) , A_j(2) \right] - \left[A_j(1) , A_i(2) \right] &= (i-j) A_{i+j-1}(1),\\
\left[p^*_i , A_j(2) \right] - \left[p^*_j , A_i(2) \right] &= (i-j) p^*_{i+j-1}
\end{aligned}
\end{equation}
which are \eqref{CommutatorAA} with $s=1, s'=2$ and \eqref{CommutatorpA} for $s=2$.

Moreover, 
\begin{equation}
\begin{aligned}
    \left[A_i(2) , A_j(2) \right] &= \left[\sum_{n\geq 1} J^{(b)}_{i-n} A_n(1) + b(i-1) A_i(1), A_j(2)\right]\\
    &\begin{multlined} = \biggl(\sum_{n\geq j} (n-j) + \sum_{n\geq i} (i-n) + \sum_{n=1}^{i+j-1} (i-n)\biggr) J^{(b)}_{i+j-n-1} A_n(1) \\ + b(i-j)(i+j-2) A_{i+j-1}(1)\end{multlined}
\end{aligned}
\end{equation}
On one hand,
\begin{multline}
\biggl(\sum_{n\geq j} (n-j) + \sum_{n\geq i} (i-n)\biggr) J^{(b)}_{i+j-n-1} A_n(1) \\
= \biggl(\sum_{n\geq M} (i-j) + \sgn(i-j) \sum_{n=\mu}^{M-1} (n-\mu)\biggr) J^{(b)}_{i+j-n-1} A_n(1).
\end{multline}
In the second sum, we notice that the sum can be started at $n=\mu+1$, we write $A_n(1) = J^{(b)}_{i-1}/(1+b)$, perform the change of summation index $n\to i+j-n$ and commute the two $J^{(b)}$s, so that
\begin{equation}
\begin{aligned}
\sum_{n=\mu+1}^{M-1} (n-\mu)J^{(b)}_{i+j-n-1} A_n(1) &= \sum_{n=\mu+1}^{M-1} (M-n) J^{(b)}_{i+j-n-1} A_n(1)\\
&= \frac{M-n}{2} \sum_{n=\mu+1}^{M-1} J^{(b)}_{i+j-n-1} A_n(1)
\end{aligned}
\end{equation}
the last equality being obtained by taking the half-sum of both expressions. This gives
\begin{multline}
\biggl(\sum_{n\geq j} (n-j) + \sum_{n\geq i} (i-n)\biggr) J^{(b)}_{i+j-n-1} A_n(1) = (i-j) \biggl(\sum_{n\geq M}  + \frac{1}{2} \sum_{n=\mu+1}^{M-1} \biggr) J^{(b)}_{i+j-n-1} A_n(1).
\end{multline}
On the other hand, using the same tricks as above, one finds
\begin{equation}
\sum_{n=1}^{i+j-1} (i-n) J^{(b)}_{i+j-n-1} A_n(1) = \frac{i-j}{2} \sum_{n=1}^{i+j-1} J^{(b)}_{i+j-n-1} A_n(1)
\end{equation}
which is decomposed as
\begin{equation}
\frac{i-j}{2} \sum_{n=1}^{i+j-1} J^{(b)}_{i+j-n-1} A_n(1) = \frac{i-j}{2} \biggl(\sum_{n=1}^{\mu} + \sum_{n=\mu+1}^{M-1} + \sum_{n=M}^{i+j-1}\biggr) J^{(b)}_{i+j-n-1} A_n(1).
\end{equation}
The first and third sum can be shown to be equal, by using $A_n(1) = J^{(b)}_{n-1}/(1+b)$, performing the change $n\to i+j-n$ and commuting the two $J^{(b)}$s. All in all, it comes that
\begin{equation}
\left[A_i(2) , A_j(2) \right] = (i-j) A_{i+j-1}(2),
\end{equation}
that is \eqref{CommutatorA} for $s=2$.

\subsection{Computation of $[ A_i(1), A_j(3) ]$} \label{sec:CommutatorA1A3}
We start with $\left[ J^{(b)}_j, A_i(3)\right]$ and use \eqref{eq:mod_rec_no_ui},
\begin{equation}
\left[ J^{(b)}_j, A_i(3)\right] = \sum_{n\geq 1} J_{i-n}^{(b)} [J^{(b)}_j, A_n(2)] + [J^{(b)}_j, J^{(b)}_{i-n}]A_n(2) + b(i-1)[J^{(b)}_j, A_i(2)].
\end{equation}
All commutators have been evaluated previously. This gives
\begin{multline}
\left[ J^{(b)}_j, A_i(3)\right] = j\left(\delta_{j<0} \sum_{l\geq j+1} + \sum_{l\geq \max(1, j+1)} + \delta_{i+j\geq 1}\sum_{l\geq 1}\right) J^{(b)}_{i+j-l} J^{(b)}_{l-1}\\
+bj\left( (i-1-j)\delta_{j<0} + (j+2i-2)\delta_{i+j\geq 1}\right) J^{(b)}_{i+j-1} - b^2 (i-1)^3 \delta_{i+j,1}.
\end{multline}
Therefore, by taking $j>0$ one finds
\begin{equation}
\left[p_j^*, A_i(3)\right] = j\left(\sum_{l\geq j} + \sum_{l\geq 1}\right) J^{(b)}_{i+j-l-1} p^*_{l} + bj(j+2i-2) p^*_{i+j-1},
\end{equation}
and therefore
\begin{multline}
\left[p_i^*, A_j(3)\right] - \left[ p_j^*, A_i(3)\right] = b(i-j)(i+j-2) p^*_{i+j-1}\\
+ \left(i\sum_{l\geq i} - j\sum_{l\geq j} + (i-j)\sum_{l\geq 1}\right) J^{(b)}_{i+j-1-l} p_l^*.
\end{multline}
Using $M=\max(i,j)$ and $\mu=\min(i, j)$ as before, it leads to
\begin{multline} \label{Explicitp*A3}
\left[p_i^*, A_j(3)\right] - \left[ p_j^*, A_i(3)\right] = b(i-j)(i+j-2) p^*_{i+j-1}\\
+ \left(2(i-j)\sum_{l\geq M} + (i-j)\sum_{l=1}^{M-1} - \operatorname{sgn}(i-j) \mu\sum_{l=\mu}^{M-1} \right) J^{(b)}_{i+j-1-l} p_l^*.
\end{multline}
as well as
\begin{multline} \label{ExplicitA1A3}
\left[ A_i(1), A_j(3)\right] - \left[ A_j(1), A_i(3)\right] = b(i-j)(i+j-2) A_{i+j-1}(1)\\
+ \left((i-j)\sum_{l\geq M} + (i-j)\sum_{l\geq 1} - \operatorname{sgn}(i-j) \sum_{l=\mu}^{M-1} (\mu-1)\right) J^{(b)}_{i+j-1-l} A_l(1).
\end{multline}

These quantities should be respectively given by the RHS of Equations \eqref{CommutatorpA} and \eqref{CommutatorAA}. Explicitly,
\begin{multline} \label{Thmp*A3}
\sum_{l\geq 1} D_{ij, l}(3) p_l^* = b(i-j)(i+j-2) p^*_{i+j-1}\\
+ \left(2(i-j)\sum_{l\geq M} + (i-j)\sum_{l=M}^{i+j-1} + \operatorname{sgn}(i-j)\sum_{l=\mu}^{M-1}(2l-3\mu+1)\right)J^{(b)}_{i+j-1-l} p_l^*
\end{multline}
and
\begin{multline} \label{ThmA1A3}
\sum_{l\geq 1} D_{ij, l}(1) A_l(3) + D_{ij, l}(3) A_l(1) = b(i-j)(i+j-2) A_{i+j-1}(1)\\
+ \left(2(i-j)\sum_{l\geq M} + (i-j)\sum_{l=M}^{i+j-1} + \operatorname{sgn}(i-j)\sum_{l=\mu}^{M-1}(2l-3\mu+1)\right)J^{(b)}_{i+j-1-l} A_l(1).
\end{multline}
Let us prove explicitly that the RHS of \eqref{ExplicitA1A3} and \eqref{ThmA1A3} are the same. Taking advantage of the fact that $A_l(1) = J^{(b)}_{l-1}/(1+b)$, we can make the change of summation index $l\to i+j-l$ and commute the two $J$s, so that
\begin{equation}
\sum_{l=\mu}^{M-1}(2l-3\mu+1)J^{(b)}_{i+j-1-l} J^{(b)}_{l-1} = \sum_{l=\mu+1}^M (2M-\mu-2l+1)J^{(b)}_{i+j-1-l} J^{(b)}_{l-1},
\end{equation}
and by then taking the half-sum,
\begin{equation}
\sum_{l=\mu}^{M-1}(2l-3\mu+1)J^{(b)}_{i+j-1-l} J^{(b)}_{l-1} = \left((M-\mu)\sum_{l=\mu+1}^{M-1} - (\mu-1)\sum_{l=\mu}^{M-1}\right)J^{(b)}_{i+j-1-l} J^{(b)}_{l-1}.
\end{equation}
The difference between the RHS of \eqref{ThmA1A3} and that of \eqref{ExplicitA1A3} is thus
\begin{multline}
\frac{(i-j)}{1+b}\left( \sum_{l\geq M} + \sum_{l=M}^{i+j-1} + \sum_{l=\mu+1}^{M-1} - \sum_{l\geq 1}\right)J^{(b)}_{i+j-1-l} J^{(b)}_{l-1}\\
= \frac{(i-j)}{1+b}\left(\sum_{l=\mu+1}^{i+j-1} - \sum_{l=1}^{M-1}\right)J^{(b)}_{i+j-1-l} J^{(b)}_{l-1}
\end{multline}
We consider the sum from $\mu+1$ to $i+j-1$ and change the summation index to $l\to i+j-l$ and commute the two $J$s which gives
\begin{equation}
\sum_{l=\mu+1}^{i+j-1} J^{(b)}_{i+j-1-l} J^{(b)}_{l-1} = \sum_{l=1}^{M-1} J^{(b)}_{i+j-1-l} J^{(b)}_{l-1},
\end{equation}
so that the above quantity vanishes and we have $\left[ A_i(1), A_j(3)\right] - \left[ A_j(1), A_i(3)\right] = \sum_{l\geq 1} D_{ij, l}(1) A_l(3) + D_{ij, l}(3) A_l(1)$, which is \eqref{CommutatorAA} for $s=1, s'=3$. The exact same steps can be followed to prove $\left[ p_i^*, A_j(3)\right] - \left[ p_j^*, A_i(3)\right] = \sum_{l\geq 1} D_{ij, l}(3) p_l^*$.

\subsection{Computation of $[A_i(2), A_j(3) ]$}
By writing again $A_j(3) = \sum_{l\geq 1} J^{(b)}_{j-l} A_l(2) + b(j-1) A_j(2)$, one finds
\begin{multline}
[ A_i(2), A_j(3) ] = \left(\sum_{n\geq j} (n-j) + \sum_{n=1}^{i+j-1} (n-j) + \sum_{n\geq i} (2i-n-1)\right) J^{(b)}_{i+j-1-n} A_n(2) \\
+ b((i-1)^2 + (j-1)(i-j)) A_{i+j-1}(2).
\end{multline}
It is then straightforward to write
\begin{multline}
[ A_i(2), A_j(3) ] - [ A_j(2), A_i(3) ] = (i-j) A_{i+j-1}(3) + b(i-j)(i+j-2)A_{i+j-1}(2) \\
+ \left(2(i-j)\sum_{n\geq M} + (i-j) \sum_{n=M}^{i+j-1} + \operatorname{sgn}(i-j)\sum_{n=\mu}^{M-1}(2n-3\mu+1)\right) J^{(b)}_{i+j-1-n} A_n(2),
\end{multline}
i.e. $[ A_i(2), A_j(3) ] - [ A_j(2), A_i(3) ] = \sum_{n\geq 1} D_{ij, n}(2) A_n(3) + D_{ij, n}(3) A_n(2)$ as desired.

\subsection{Computation of $[ A_i(3), A_j(3) ]$} \label{sec:A3A3}
The last element of the proof of Proposition \ref{thm:SimplifiedCommutators} is the computation of $[ A_i(3), A_j(3) ]$, which is quite intricate. We provide a rather detailed proof in the present section. For bookkeeping purposes, we will use a notation which keeps track of both the origins and the types of all different contributions entering the calculation. The origin of a term will be denoted by a letter, $B, C, D, E, F, G, H$. The type of a term will refer to its the form of its summand. An additional number may appear if there are several terms with the same origin and of the same type. For instance,
\begin{equation*}
C^{JJA,2}_{ij} \coloneqq -(i-j)\sum_{\ell =1}^{M-\mu} \sum_{n=\mu+\ell}^{M-1} J^{(b)}_{i+j-1-n} J^{(b)}_{n-\ell} A_\ell(2)
\end{equation*}
has an origin indicated by the letter $C$, its type is $JJA$ (refering to its summand) and it is the second term of these origins and types.

We split the commutator as follows
\begin{multline}
    \left[A_i(3), A_j(3) \right] = \sum\limits_{n,\ell \geq 1} \underbrace{ \left[J^{(b)}_{i-n}A_n(2), J^{(b)}_{j-\ell} A_\ell(2)\right]}_{Q_{ij}} + \underbrace{b^2(i-1)(j-1) \left[A_i(2), A_j(2) \right]}_{H_{ij}} \\
    + \underbrace{b(i-1)\sum\limits_{\ell \geq 1}  \left[A_i(2), J^{(b)}_{j-\ell}A_\ell(2) \right]}_{G_{ij}} + \underbrace{b(j-1)\sum\limits_{n \geq 1}  \left[J^{(b)}_{i-n}A_n(2), A_j(2) \right]}_{-G_{ji}}
\end{multline}
We further expand the term $Q_{ij}$ as
\begin{multline}
    Q_{ij} = \underbrace{\sum\limits_{n,\ell \geq 1} J^{(b)}_{i-n}J^{(b)}_{j-\ell}\left[A_n(2), A_\ell(2)\right] }_{C_{ij}} + \underbrace{\sum\limits_{n,\ell \geq 1} \left[J^{(b)}_{i-n},J^{(b)}_{j-\ell}\right] A_\ell(2) A_n(2)}_{B_{ij}} \\
    + \underbrace{\sum\limits_{n,\ell \geq 1} J^{(b)}_{i-n}\left[A_n(2),J^{(b)}_{j-\ell}\right] A_\ell(2) }_{P_{ij}} + \underbrace{\sum\limits_{n,\ell \geq 1} J^{(b)}_{j-\ell}\left[J^{(b)}_{i-n}, A_\ell(2)\right]A_n(2)  }_{-P_{ji}}.
\end{multline}
We now analyze all terms so as to ``reform'' $\sum_{n\geq 1} D_{ij,n}(3) A_n(3)$.

\subsubsection{Term $B_{ij}$}
It writes
\begin{equation}
    B_{ij} = (1+b)\sum\limits_{n=1}^{i+j-1} (i-n) A_{i+j-n}(2) A_n(2).
\end{equation}
We make the change of summation index $n\to i+j-n$ and commute the two $A(2)$,
\begin{equation}
B_{ij} = (1+b) \sum\limits_{n=1}^{i+j-1} (n-j) A_{i+j-n}(2) A_n(2) + \frac{1+b}{6} (i+j)(i+j-1)(i+j-2) A_{i+j-1}(2).
\end{equation}
We then take the half-sum
\begin{equation}
     B_{ij} = \underbrace{\frac{1+b}{2} (i-j)\sum\limits_{n=1}^{i+j-1} A_{i+j-n}(2) A_n(2)}_{B^{AA}_{ij}} + \underbrace{\frac{1+b}{12} (i+j)(i+j-1)(i+j-2) A_{i+j-1}(2)}_{B^{A}_{ij}}
\end{equation}

\subsubsection{Term $C_{ij}$} \label{sec:C}
It expands as
\begin{equation}
C_{ij} = \sum_{n, l\geq 1} (n-l) J^{(b)}_{i-n} J^{(b)}_{j-l} A_{n+l-1}(2) = \sum_{k\geq 1} \sum_{n=1}^k (2n-k-1) J^{(b)}_{i-n} J^{(b)}_{j+n-k-1} A_k(2).
\end{equation}
We can either make the change $p=j+n-1$ and get
\begin{equation}
C_{ij} = \sum_{k\geq 1} \sum_{p=j}^{j+k-1} (2p-2j-k+1) J^{(b)}_{i+j-p-1} J^{(b)}_{p-k} A_k(2)
\end{equation}
or $p=i+k-n$ and get
\begin{equation}
C_{ij} = \sum_{k\geq 1} \sum_{p=i}^{i+k-1} (2i+k-2p-1)  J^{(b)}_{p-k} J^{(b)}_{i+j-p-1} A_k(2)
\end{equation}
By commuting the two $J$s the latter gives
\begin{equation}
C_{ij} = \sum_{k\geq 1} \sum_{p=i}^{i+k-1} (2i+k-2p-1) J^{(b)}_{i+j-p-1} J^{(b)}_{p-k} A_k(2) \underbrace{- \frac{1+b}{6} (i+j)(i+j-1)(i+j-2) A_{i+j-1}(2)}_{2C_{ij}^{A,1}}.
\end{equation}
We then take the half-sum and split the terms as follows
{\small
\begin{multline}
C_{ij} = C_{ij}^{A,1} + \underbrace{ (i-j) \sum_{\ell \geq M-\mu+1} \sum_{n=M}^{\ell+\mu-1} J^{(b)}_{i+j-1-n} J^{(b)}_{n-\ell} A_\ell(2)}_{C^{JJA,1}_{ij}} \\
+ \frac{\operatorname{sgn}(j-i)}{2} \left[\underbrace{\sum_{\ell \geq 1} \sum_{n=\mu}^{M-1}}_{S_1} \underbrace{- \sum_{\ell =1}^{M-\mu} \sum_{n=\mu+\ell}^{M-1} }_{T_1}\right] (2\mu-2n-1+\ell) J^{(b)}_{i+j-1-n} J^{(b)}_{n-\ell} A_\ell(2) \\
- \frac{\operatorname{sgn}(j-i)}{2} \left[\underbrace{\sum_{\ell \geq 1} \sum_{n=\mu+\ell}^{M+\ell-1}}_{S_2} \underbrace{- \sum_{\ell =1}^{M-\mu} \sum_{n=\mu+\ell }^{M-1}}_{T_2}\right] (2M-2n-1+\ell) J^{(b)}_{i+j-1-n} J^{(b)}_{n-\ell} A_\ell(2)
\end{multline}}%
where we can already notice $C_{ij}^{A,1} = -B_{ij}^{A}$.

The sums $T_1+T_2$ get together as
\begin{equation}
C^{JJA,2}_{ij} \coloneqq -(i-j)\sum_{\ell =1}^{M-\mu} \sum_{n=\mu+\ell}^{M-1} J^{(b)}_{i+j-1-n} J^{(b)}_{n-\ell} A_\ell(2).
\end{equation}
In order to pack together $S_1$ and $S_2$, we relabel $n \rightarrow i+j-1-n+\ell$ in $S_2$ and commute the two $J$s. It comes
\begin{multline}
C_{ij} = C^{JJA,1}_{ij} + C^{JJA,2}_{ij} + C_{ij}^{A,1} + \underbrace{\frac{\operatorname{sgn}(j-i)}{2}\sum_{n=\mu}^{M-1} (2n-3\mu-M+2)(i+j-1-n) A_{i+j-1}(2)}_{C^{A,2}_{ij}} \\
 \underbrace{-\operatorname{sgn}(j-i) \sum_{\ell \geq 1} \sum_{n=\mu}^{M-1} (2n-2\mu+1-\ell) J^{(b)}_{i+j-1-n} J^{(b)}_{n-\ell} A_\ell(2)}_{\tilde{C}^{JJA}_{ij}}.
\end{multline}

\subsubsection{Term $P_{ij}$} \label{sec:P}
This contribution splits into three terms,
\begin{multline}
P_{ij} = \underbrace{\sum\limits_{n \geq j} b(n-j)^2J^{(b)}_{i+j-1-n} A_{n}(2)}_{F_{ij}} 
+ \underbrace{\sum\limits_{n, \ell \geq j}  (\ell-j) J^{(b)}_{i+j-1-n} J^{(b)}_{n-\ell} A_\ell(2)}_{D_{ij}} \\
+ \underbrace{\sum\limits_{n \geq j} \sum_{\ell=1}^{n} (\ell-j) J^{(b)}_{i+j-1-n} J^{(b)}_{n-\ell} A_\ell(2)}_{E_{ij}}.
\end{multline}
Recall that we do not need $P_{ij}$ itself but $P_{ij}-P_{ji}$, so we can compute directly the antisymmetrized contribution for each of the terms above. We have
\begin{multline}
    D_{ij}-D_{ji} = (i-j)\sum\limits_{n, \ell \geq M} J^{(b)}_{i+j-1-n} J^{(b)}_{n-\ell}A_\ell(2) \\
    +  \underbrace{\sgn(j-i) \left[ \sum\limits_{n \geq M}  \sum\limits_{\ell = \mu}^{M-1} +  \sum\limits_{n = \mu}^{M-1} \sum\limits_{\ell = \mu}^{M-1} +  \sum\limits_{n = \mu}^{M-1}  \sum\limits_{\ell \geq M}\right] (\mu-\ell) J^{(b)}_{i+j-1-n} J^{(b)}_{n-\ell} A_\ell(2)}_{\Omega^{(1)}_{ij}}.
\end{multline}
We use the first term above to form some $A_{n}(3)$,
\begin{multline}
(i-j)\sum\limits_{n, \ell \geq M} J^{(b)}_{i+j-1-n} J^{(b)}_{n-\ell}A_\ell(2) = (i-j)\sum\limits_{n \geq M}  J^{(b)}_{i+j-1-n} A_n(3) \\
\underbrace{-(i-j) \sum\limits_{n \geq M} \sum\limits_{\ell =1}^{M-1} J^{(b)}_{i+j-1-n} J^{(b)}_{n-\ell} A_\ell(2)}_{\alpha_{ij}}
\underbrace{-b(i-j)\sum\limits_{n \geq M} (n-1) J^{(b)}_{i+j-n-1} A_n(2)}_{D^{bJA,1}_{ij}}
\end{multline}

Similarly the sum over $n \geq M$ in $\alpha_{ij}$ can be completed to give a contribution of the form $A(2)A(2)$,
\begin{multline}
    \alpha_{ij} = - (1+b)(i-j) \sum_{\ell=1}^{M-1} A_{i+j-\ell}(2) A_\ell(2) + \underbrace{(i-j) \sum_{\ell=1}^{M-1} \sum_{n = \ell}^{M-1} J^{(b)}_{i+j-n-1} J^{(b)}_{n-\ell} A_\ell(2)}_{D^{JJA}_{ij}} \\
    + \underbrace{b(i-j) \sum_{\ell=1}^{M-1} (i+j-1-\ell) J^{(b)}_{i+j-1-\ell} A_\ell(2)}_{D^{bJA,2}_{ij}}
\end{multline}
In the first term, we commute the two $A(2)$, relabel $\ell \rightarrow i+j-\ell$ and take the half-sum, leading to 
\begin{multline}
    \alpha_{ij} = D^{JJA,1}_{ij} + D^{bJA,2}_{ij} \underbrace{-\frac{1+b}{2} (i-j) \sum_{\ell=1}^{i+j-1}A_{i+j-\ell}(2) A_\ell(2)}_{D^{AA,1}_{ij}} \\
     \underbrace{-\frac{1+b}{2} (i-j) \sum_{\ell=\mu+1}^{M-1} A_{i+j-\ell}(2) A_\ell(2)}_{D^{AA,2}_{ij}} + \underbrace{ \frac{1+b}{2}(i-j) \sum_{\ell=1}^{M-1} (2\ell-i-j) A_{i+j-1}(2)}_{D^A_{ij}}
\end{multline}

Finally, we split $\Omega_{ij}^{(1)}$ as 
\begin{multline}
\Omega_{ij}^{(1)} = \underbrace{\sgn(j-i) \sum\limits_{n \geq M}  \sum\limits_{\ell = \mu}^{M-1} (\mu-\ell) J^{(b)}_{i+j-1-n} J^{(b)}_{n-\ell} A_\ell(2)}_{\tilde{D}^{JJA,2}_{ij}} \\
+ \underbrace{\sgn(j-i) \sum\limits_{n = \mu}^{M-1} \sum\limits_{\ell \geq \mu} (\mu-\ell) J^{(b)}_{i+j-1-n} J^{(b)}_{n-\ell} A_\ell(2)}_{\tilde{D}^{JJA,1}_{ij}}
\end{multline}
All in all, we find the decomposition
\begin{multline}
D_{ij}-D_{ji}
= (i-j)\sum\limits_{n \geq M}  J^{(b)}_{i+j-1-n} A_n(3) \\
+ D^{AA,1}_{ij} + D^{AA,2}_{ij} + D^A_{ij} + D^{JJA}_{ij} + \tilde{D}^{JJA,1}_{ij} + \tilde{D}^{JJA,2}_{ij} + D^{bJA,1}_{ij} + D^{bJA,2}_{ij}
\end{multline}

For the term $E_{ij}$ we get
\begin{multline}
    E_{ij} - E_{ji} = \underbrace{(i-j)\sum\limits_{n \geq M} \sum\limits_{\ell = 1}^{n-1} J^{(b)}_{i+j-1-n}J^{(b)}_{n-\ell} A_\ell(2) }_{X_{ij}}\\
    + \underbrace{\sgn(j-i) \sum\limits_{n =\mu}^{M-1}\sum\limits_{\ell=1}^{n-1} (\mu-\ell) J^{(b)}_{i+j-1-n}J^{(b)}_{n-\ell} A_\ell(2)}_{\Omega^{(2)}_{ij}}
\end{multline}
The latter term is split as follows
\begin{multline}
\Omega^{(2)}_{ij} = \underbrace{\sgn(j-i) \sum\limits_{n =\mu}^{M-1}\sum\limits_{\ell=1}^{\mu-1} (\mu-\ell) J^{(b)}_{i+j-1-n}J^{(b)}_{n-\ell} A_\ell(2)}_{\tilde{E}^{JJA,1}_{ij}} \\
+ \underbrace{\sgn(j-i) \sum\limits_{n =\mu}^{M-1}\sum\limits_{\ell=\mu}^{n} (\mu-\ell) J^{(b)}_{i+j-1-n}J^{(b)}_{n-\ell} A_\ell(2)}_{\tilde{E}^{JJA,2}_{ij}}
\end{multline}
We complete the sum on $\ell$ to get a $A_n(3)$,
\begin{multline}
X_{ij} = (i-j) \sum\limits_{n \geq M} J^{(b)}_{i+j-1-n} A_{n}(3) \underbrace{- (i-j) \sum\limits_{n \geq M} \sum\limits_{\ell \geq n} J^{(b)}_{i+j-1-n} J^{(b)}_{n-\ell} A_{\ell}(2)}_{E^{JJA}_{ij}} \\
\underbrace{-b(i-j) \sum\limits_{n \geq M} (n-1) J^{(b)}_{i+j-1-n} A_{n}(2)}_{E^{bJA}_{ij}}
\end{multline}
hence
\begin{equation}
E_{ij} - E_{ji} = (i-j) \sum\limits_{n \geq M} J^{(b)}_{i+j-1-n} A_{n}(3) + \tilde{E}^{JJA,1}_{ij} + \tilde{E}^{JJA,2}_{ij} + E^{JJA}_{ij} + E^{bJA}_{ij}.
\end{equation}

And $F_{ij}$ gives
\begin{multline}
    F_{ij} - F_{ji} = \underbrace{b\sum\limits_{n \geq M} (i-j)(2n-i-j) J^{(b)}_{i+j-1-n} A_{n}(2)}_{F^{bJA}_{ij}} \\
    + \underbrace{b\sgn(i-j) \sum\limits_{n=\mu}^{M-1} (n-\mu)^2 J^{(b)}_{i+j-1-n} A_{n}(2)}_{\tilde{F}^{bJA}_{ij}}.
\end{multline}

\subsubsection{Terms $G_{ij}$ and $H_{ij}$}
The term $H_{ij} = b^2(i-1)(j-1) \left[A_i(2), A_j(2) \right]$ is simply
\begin{equation}
    H_{ij} = b^2(i-1)(j-1)(i-j) A_{i+j-1}(2).
\end{equation}

As for $G_{ij} = b(i-1)\sum_{l\geq 1} [A_i(2), J^{(b)}_{j-l} A_l(2)]$,
\begin{multline}
    G_{ij} = \frac{b^2}{1+b}(i-1)^3 A_{i+j-1}(2) \\
    + b(i-1) \left(\sum_{l\geq i} (2i-l-1) + \sum_{l\geq j} (l-j) + \sum_{l=1}^{i+j-1} (l-j)\right) J^{(b)}_{i+j-l-} A_l(2).
\end{multline}
Antisymmetrization leads to
\begin{multline}
    G_{ij}-G_{ji} = \underbrace{b^2\left((i-1)^3-(j-1)^3\right) A_{i+j-1}(2)}_{G^{A}_{ij}} \\
    + \underbrace{b(i-j)\sum\limits_{n=1}^{i+j-1} (n-1) J^{(b)}_{i+j-1-n} A_{n}(2)}_{G^{bJA,1}_{ij}}
    + \underbrace{2b(i-j)(i+j-2)\sum\limits_{n \geq M}  J^{(b)}_{i+j-1-n} A_{n}(2)}_{G^{bJA,2}_{ij}} \\
    + \underbrace{b\sgn(i-j) \sum\limits_{n=\mu}^{M-1} \left[(M-1)(n-\mu) - (\mu-1)(2\mu-n-1)\right] J^{(b)}_{i+j-1-n} A_{n}(2)}_{\tilde{G}^{bJA}_{ij}}
\end{multline}

\subsubsection{Repackaging the contributions} \label{sec:Repackage}
We start with adding together the three terms $\tilde{C}^{JJA}_{ij}, \tilde{D}^{JJA,1}_{ij}, \tilde{E}^{JJA,1}_{ij}$.
\begin{equation}
\begin{aligned}
&\tilde{C}^{JJA}_{ij}+ \tilde{D}^{JJA,1}_{ij}+ \tilde{E}^{JJA,1}_{ij}
= \sgn(j-i)\\
& \sum_{n=\mu}^{M-1} \left(\sum_{l\geq \mu} (\mu-l) + \sum_{l=1}^{\mu-1} (\mu-l) - \sum_{l\geq 1}(2n-2\mu-l+1)\right)J^{(b)}_{i+j-1-n} J^{(b)}_{n-l} A_l(2)\\
&= \sgn(j-i) \sum_{n=\mu}^{M-1} \sum_{l\geq 1} (3\mu-2n-1) J^{(b)}_{i+j-1-n} J^{(b)}_{n-l} A_l(2)\\
&= \sgn(j-i) \sum_{n=\mu}^{M-1} (3\mu-2n-1) J^{(b)}_{i+j-1-n} (A_n(3) - b (n-1) A_n(2))\\
&= \sgn(j-i) \sum_{n=\mu}^{M-1} (3\mu-2n-1) J^{(b)}_{i+j-1-n} A_n(3) + \widetilde{{CE}}^{bJA}_{ij}
\end{aligned}
\end{equation}
with $\widetilde{{CE}}^{bJA}_{ij} = - b\sgn(j-i) \sum_{n=\mu}^{M-1} (3\mu-2n-1)(n-1) J^{(b)}_{i+j-1-n} A_n(2)$. Then
\begin{equation}
\begin{aligned}
\tilde{D}^{JJA,2}_{ij} + \tilde{E}^{JJA,2}_{ij} &= \sgn(j-i) \sum_{l=\mu}^{M-1}(\mu-l) \sum_{n\geq l} J^{(b)}_{i+j-1-n} J^{(b)}_{n-l} A_l(2)\\
&= \sgn(j-i) \sum_{l=\mu}^{M-1}(\mu-l) \sum_{p\geq 1} J^{(b)}_{i+j-1-p-l} J^{(b)}_{p-1} A_l(2)\\
&= \sgn(j-i) (1+b) \sum_{l=\mu}^{M-1}(\mu-l) (A_{i+j-l}(2) - b (i+j-l-1) A_{i+j-l}(1)) A_l(2).
\end{aligned}
\end{equation}
The term with $A_{i+j-l}(2) A_l(2)$ can be rewritten by making the change of index variable $l\to i+j-l$ and commuting the two operators,
\begin{multline}
\sgn(j-i) (1+b) \sum_{l=\mu+1}^{M-1}(\mu-l) A_{i+j-l}(2) A_l(2)
= \sgn(j-i) (1+b) \sum_{l=\mu+1}^{M-1}(l-M) A_{i+j-l}(2) A_l(2) \\+ \sgn(j-i) (1+b) \sum_{l=\mu+1}^{M-1}(l-M)(2l-i-j) A_{i+j-1}(2).
\end{multline}
then by taking the half-sum, we get
\begin{multline}
\tilde{D}^{JJA,2}_{ij} + \tilde{E}^{JJA,2}_{ij} = \underbrace{\frac{1+b}{2} (i-j) \sum_{l=\mu+1}^{M-1} A_{i+j-l}(2) A_l(2)}_{(DE)^{AA}_{ij}}\\
+ \underbrace{\sgn(i-j) b\sum_{l=\mu}^{M-1} (\mu-l)(i+j-l-1) J^{(b)}_{i+j-l-1} A_l(2)}_{\widetilde{(DE)}^{bJA}_{ij}}\\
+ \underbrace{\sgn(j-i) \frac{1+b}{2} \sum_{l=\mu+1}^{M-1}(l-M)(2l-i-j) A_{i+j-1}(2)}_{(DE)^A_{ij}}
\end{multline}

Next up, we look at $C^{JJA,1}_{ij} + E^{JJA}_{ij}$,
\begin{equation}
\begin{aligned}
C^{JJA,1}_{ij} + E^{JJA}_{ij} &= (i-j) \sum_{n\geq M} \sum_{l=n-\mu+1}^{n-1} J^{(b)}_{i+j-1-n} J^{(b)}_{n-l} A_{l}(2)\\
&= (i-j) \Biggl( \underbrace{\sum_{l\geq M-1} \sum_{n=l+1}^{l+\mu-1}}_{z_{ij}} + \underbrace{\sum_{l=M-\mu+1}^{M-2} \sum_{n=M}^{l+\mu-1}}_{(CE)^{JJA,1}_{ij}}\Biggr) J^{(b)}_{i+j-1-n} J^{(b)}_{n-l} A_{l}(2)
\end{aligned}
\end{equation}
then
\begin{equation}
\begin{aligned}
z_{ij} &= (i-j) \sum_{l\geq M-1} \sum_{p=1}^{\mu-1} J^{(b)}_{i+j-1-l-p} J^{(b)}_{p} A_{l}(2)\\
&= \underbrace{(i-j) \sum_{l\geq M-1} \sum_{p=1}^{\mu-1}  J^{(b)}_{p} J^{(b)}_{i+j-1-l-p} A_{l}(2)}_{\zeta_{ij}} \underbrace{- (i-j) (1+b)\frac{\mu(\mu-1)}{2} A_{i+j-1}(2)}_{(CE)^A_{ij}}
\end{aligned}
\end{equation}
In the first line we have set $p=n-l$, in the second line we have commuted the two $J$s. In $\zeta_{ij}$ we form $A_{i+j-1-p}(3)$
\begin{multline}
\zeta_{ij} = (i-j) \sum_{p=1}^{\mu-1} J^{(b)}_{p} (A_{i+j-p-1}(3) - b(i+j-p-2) A_{i+j-1-p}(2))\\
- (i-j) \sum_{p=1}^{\mu-1} \sum_{l=1}^{M-2} J^{(b)}_{p} J^{(b)}_{i+j-1-l-p} A_{l}(2)
\end{multline}
and making the change $n=i+j-1-p$,
\begin{multline}
\zeta_{ij} = (i-j) \sum_{n=M}^{i+j-1} J^{(b)}_{i+j-1-n} A_{n}(3) \underbrace{- b(i-j) \sum_{n=M}^{i+j-1} (n-1) J^{(b)}_{i+j-1-n} A_{n}(2)}_{(CE)^{bJA}_{ij}}\\
\underbrace{- (i-j) \sum_{n=M}^{i+j-2} \sum_{l=1}^{M-2} J^{(b)}_{i+j-1-n} J^{(b)}_{n-l} A_{l}(2)}_{(CE)^{JJA,2}_{ij}}.
\end{multline}
Overall
\begin{multline}
C^{JJA,1}_{ij} + E^{JJA}_{ij} = (i-j) \sum_{n=M}^{i+j-1} J^{(b)}_{i+j-1-n} A_{n}(3) \\
+ (CE)^{JJA,1}_{ij} + (CE)^{JJA,2}_{ij} + (CE)^A_{ij} + (CE)^{bJA}_{ij}.
\end{multline}


We package together some of the contributions which receive some $b$ during the calculation. First,
\begin{equation}
    H_{ij} + G^{A}_{ij} = b^2(i-j)(i+j-2)^2 A_{i+j-1}(2).
\end{equation}
Moreover,
\begin{align}
D^{bJA,1}_{ij} + E^{bJA}_{ij} + F^{bJA}_{ij} + G^{bJA,2}_{ij} &= b(i-j)(i+j-2) \sum_{n\geq M} J^{(b)}_{i+j-1-n} A_n(2)\\
(CE)^{bJA}_{ij} + G^{bJA,1}_{ij} + D^{bJA,2}_{ij} &= b(i-j)(i+j-2) \sum_{n=1}^{M-1} J^{(b)}_{i+j-1-n} A_n(2)
\end{align}
hence the sum of the above three contributions reduces to
\begin{multline}
(D^{bJA,1}_{ij} + E^{bJA}_{ij} + F^{bJA}_{ij} + G^{bJA,2}_{ij}) + ((CE)^{bJA}_{ij} + G^{bJA,1}_{ij} + D^{bJA,2}_{ij}) + (H_{ij} + G^{A}_{ij}) \\
= b(i-j)(i+j-2) A_{i+j-1}(3).
\end{multline}

\subsubsection{Cancellations}
We notice the following cancellations,
\begin{equation}
\begin{aligned}
C^{A,1}_{ij} + B^{A}_{ij} &= 0, \\
D^{AA,1}_{ij} + B^{AA}_{ij} &= 0, \\
D^{AA,2}_{ij} + (DE)^{AA}_{ij} &=0, \\
C^{A,2}_{ij} + D^A_{ij} + (DE)^A_{ij} + (CE)^A_{ij}  &= 0, \\
C^{JJA,2}_{ij} + D^{JJA}_{ij} + (CE)^{JJA,1}_{ij} + (CE)^{JJA,2}_{ij} &= 0, \\
\widetilde{(CE)}^{bJA}_{ij} + \widetilde{(DE)}^{bJA}_{ij} + \tilde{G}^{bJA}_{ij} + \tilde{F}^{bJA}_{ij}  &= 0.
\end{aligned}
\end{equation}

\subsubsection{Final expression}
Bringing together the remaining contributions, we get
\begin{multline}
\left[ A_i(3) , A_j(3) \right] = 2(i-j) \sum\limits_{n \geq M} J^{(b)}_{i+j-1-n} A_n(3)
+ (i-j) \sum\limits_{n=M}^{i+j-1} J^{(b)}_{i+j-1-n} A_n(3) \\
+ \sgn(i-j) \sum\limits_{n=\mu}^{M-1}(2n-3\mu+1) J^{(b)}_{i+j-1-n} A_n(3) + b(i-j)(i+j-2)A_{i+j-1}(3)
\end{multline}
which concludes the proof of Proposition~\ref{thm:SimplifiedCommutators}.

\section{Computation of the commutation relations for bipartite maps with black vertices of degrees at most 3}
\label{sec:commk1}
\subsection{Proof of Theorem~\ref{thm:2ndModel}}
In this section we consider the function $\tau_{k,r}$ with $k=1, r=3$, characterized by the evolution equation
\begin{equation*}
\frac{\partial \tau_{1,3}(t, \pp, \qq, u)}{\partial t} = \sum_{m=1}^3 t^{m-1} M^{(1,m)} \tau_{1,3}(t, \pp, \qq, u),
\end{equation*}
where $M^{(k, m)}$ is given in \eqref{Mkm}. Combinatorially, it is the generating series of $b$-deformed bipartite maps whose vertex of color 0 (or, say, black vertices) of degree $m\geq 1$ are weighted with $q_m$ vanishing for $m>3$, and whose vertices of color 1 (or, say, white vertices) are each counted with weight $u$. The candidate constraints, appearing in Theorem \ref{thm:2ndModel}, are
\begin{equation*}
L^{\text{bip}\leq 3}_i \coloneqq -p_i^* + \sum_{m=1}^3 q_m t^m [y_i] \Bigl(Y_+(\Lambda_Y + u)\Bigr)^m \frac{y_0}{1+b}
\end{equation*}

For this model, a simplification appears compared to the case of 2- and 3-constellations of the previous section: since there is a single variable $u$, there is in fact no need to expand $L^{\text{bip}\leq 3}_i$ on elementary symmetric functions of $u_1, \dotsc, u_k$. Instead, it is enough to set $J^{(b)}_0 = \alpha = u$ in \eqref{eq:currents}. {\bf In the rest of this section we thus set $J^{(b)}_0 = \alpha = u$.} Then the operators $M^{(1,m)}_i$ as given in \eqref{M1m} satisfy the following recursion.
\begin{proposition}
One has $M^{(1,1)}_i = \frac{J^{(b)}_{i-1}}{1+b}$ for all $i\geq 1$. Then for all $i\geq1, m\geq 1$
\begin{equation} \label{eq:mod_rec_no_ui'}
    M^{(1,m+1)}_i = \sum\limits_{n \geq 1} J^{(b)}_{i-n-1} M^{(1,m)}_n + b(i-1) M^{(1,m)}_{i-1}.
\end{equation}
(using $M^{(1,m)}_0 = 0$ for all $m\geq 1$ when necessary).
\end{proposition}

\begin{proof}
Let $Y_i$, $i\in\mathbb{Z}$, be the operator defined as $Y_i y_j = y_{i+j} \delta_{i+j\geq 0}$. Then,
\begin{equation} \label{Lambda+u}
\Lambda_Y+u = \sum_{i\in\mathbb{Z}} Y_i J^{(b)}_i + b \sum\limits_{i \geq 0} iy_i [y_i].
\end{equation}
By definition, $M^{(1,m)}_i = [y_i] (Y_+(\Lambda_y+u))^m \frac{y_0}{1+b}$ for all $i, m\geq 0$, which gives directly $M^{(1,1)}_i = \frac{J^{(b)}_{i-1}}{1+b}$ for all $i\geq 0$, and $M^{(1,m)}_0 = \frac{\delta_{m,0}}{1+b}$ for all $m\geq 0$. Then for $m\geq 1$
\begin{equation}
M^{(1,m+1)}_i = \sum_{n\geq 0} [y_i] Y_+(\Lambda_Y+u) y_n [y_n] (Y_+(\Lambda_y+u))^m \frac{y_0}{1+b} = \sum_{n\geq m} [y_i] Y_+(\Lambda_Y+u) y_n M^{(1,m)}_n
\end{equation}
and we conclude with \eqref{Lambda+u}.
\end{proof}

The structure of the calculations is completely similar to what was done in the previous section. The analog of Lemma \ref{lemma:SimplifiedCommutators} is as follows.
\begin{lemma} \label{lemma:SimplifiedCommutators2}
Assume that there exist differential operators $(\tilde{D}_{ij, l}(m))_{i,j,l, m\geq 1}$ in the variables $p_i$s such that
\begin{equation} \label{Commutatorp*M}
[p_i^*, M^{(1,m)}_j] - [p_j^*, M^{(1,m)}_j] = \sum_{l\geq 1} \tilde{D}_{ij, l}(m) p_l^*
\end{equation}
and
\begin{equation} \label{eq:SimplifiedCommutators2}
\begin{aligned}
    \left[M^{(1,m)}_i, M^{(1,m)}_j\right] &= \sum_{l\geq 1} \tilde{D}_{ij, l}(m) M^{(1,m)}_l, \\
    \left[M^{(1,m)}_i, M^{(1,m')}_j\right] -  \left[M^{(1,m)}_j, M^{(1,m')}_i\right] &= \sum_{l\geq 1} \tilde{D}_{ij, l}(m') M^{(1,m)}_l + \tilde{D}_{ij, l}(m) M^{(1,m')}_l.
\end{aligned}
\end{equation}
Then the operators $L^{\text{bip}\leq r}_i = -p_i^* + \sum_{m=1}^r q_m t^m M^{(1,m)}_i$ satisfy
\begin{equation}
\label{eq:CommutConstr2}
\left[L^{\text{bip}\leq r}_i, L^{\text{bip}\leq r}_j\right] = t \sum_{l\geq 1} \tilde{D}_{ij, l}^{(r)} L_l^{\text{bip}\leq r},
\end{equation}
with
\begin{equation}
\tilde{D}_{ij, l}^{(r)} = \sum_{m=1}^r q_m t^{m-1} \tilde{D}_{ij, l}(m).
\end{equation}
\end{lemma}

\begin{proof}
    The proof is identical to the one of Lemma~\ref{lemma:SimplifiedCommutators}.
\end{proof}

The analog of Proposition \ref{thm:SimplifiedCommutators} is then:
\begin{proposition} \label{thm:SimplifiedCommutators2}
Let $\mu = \min(i,j)$ and $M = \max(i,j)$. For $1\leq m, m'\leq 3$,
\begin{subequations} \label{SimplifiedCommutators2}
\begin{alignat}{1}
    \left[M^{(1,m)}_i, M^{(1,m)}_j\right] &= \sum_{l\geq 1} \tilde{D}_{ij, l}(m) M^{(1,m)}_l,\label{CommutatorM} \\
    \left[M^{(1,m)}_i, M^{(1,m')}_j\right] -  \left[M^{(1,m')}_i, M^{(1,m)}_j\right] &= \sum_{l\geq 1} \tilde{D}_{ij,l}(m') M^{(1,m)}_l + \tilde{D}_{ij,l}(m) M^{(1,m')}_l, \label{CommutatorMM}\\
    \left[p_i^*, M^{(1,m)}_j\right] -  \left[p_j^*, M^{(1,m)}_i\right] &= \sum_{l\geq 1} \tilde{D}_{ij,l}(m) p_l^* \label{CommutatorpM}
\end{alignat}
\end{subequations}
with $\tilde{D}_{ij, l}(1) = 0$ and $\tilde{D}_{ij, l}(2) = (i-j) \delta_{l,i+j-2}$ and
\begin{multline}
\tilde{D}_{ij, l}(3) = (i-j) (2\delta_{l\geq M-1} + \delta_{M-1\leq l\leq i+j-3}) J^{(b)}_{i+j-3-l} \\
+ \operatorname{sgn}(i-j)(2l-3\mu+3) \delta_{\mu-1\leq l\leq M-2} J^{(b)}_{i+j-3-l}
+ b(i-j)(i+j-3)\delta_{l, i+j-3}.
\end{multline}
\end{proposition}

This section contains the proof of Proposition~\ref{thm:SimplifiedCommutators2}. We compute all the commutators for all operators $M^{(1,m)}_i$ up to $m=3$. Its structure (and the results themselves) are analogous to those of Section~\ref{sec:3Constellations} for $A_i(p)$. However, it was not possible for us to prove \ref{thm:SimplifiedCommutators2} from the calculations of the previous section, nor to unify them.

In all this section, we set $J_0 = \alpha =u$.

\subsection{Computation of $\left[M_i^{(1,1)}, M_j^{(1,1)}\right]$}
We have $M^{(1,1)}_i = \frac{J^{(b)}_{i-1}}{1+b}$, then $\left[J^{(b)}_j, M_i^{(1,1)} \right] = j\delta_{i+j,1}$. Therefore $[M_i^{(1,1)}, M_j^{(1,1)}] = [M_i^{(1,1)}, p_j^*] = 0 $, which are indeed compatible with $\tilde{D}_{ij,l}(1)=0$.

\subsection{Computation of $\left[ M_i^{(1,m)}, M_j^{(1,m')}\right]$ for $m, m'\leq 2$}
We use relation~\eqref{eq:mod_rec_no_ui'} to compute the various commutators. We have
\begin{equation}
    \left[ J^{(b)}_j, M_i^{(1,2)}  \right] = (1+b) j M^{(1,1)}_{i+j-1}\left( \delta_{i+j \geq 2} +  \delta_{j \leq 0} \right) - b (i-1)(i-2) \delta_{i+j,2}.
\end{equation}
This $[p_j^*, M^{(1,2)}_i] = jJ^{(b)}_{i+j-2}$ for $i, j\geq 1$ and $[M^{(1,1)}_j, M^{(1,2)}_i] = (j-1) M^{(1,1)}_{i+j-2}$, hence
\begin{equation}
\begin{aligned}
    [M^{(1,2)}_i, p_j^*] - [M^{(1,2)}_j, p_i^*] &= (i-j) p_{i+j-2}^*\\ 
    [M^{(1,2)}_i, M^{(1,1)}_j] - [M^{(1,2)}_j,M^{(1,1)}_i] &= (i-j) M^{(1,1)}_{i+j-2}
\end{aligned}
\end{equation}
which correspond to \eqref{CommutatorMM} and \eqref{CommutatorpM} with $m'=1, m=2$. Furthermore, the calculation of $[M^{(1,2)}_i, M^{(1,2)}_j]$ follows exactly the same steps as $[A_i(2), A_j(2)]$ in Section \ref{sec:A2A2} and again produces a half-Virasoro algebra,
\begin{equation}
[M^{(1,2)}_i, M^{(1,2)}_j] = (i-j) M^{(1,2)}_{i+j-2}.
\end{equation}
which is also compatible with $\tilde{D}_{ij,l}(2) = (i-j) \delta_{l,i+j-2}$.

\subsection{Computation of $\left[ M^{(1,1)}_i, M^{(1,3)}_j\right]$}
Again, we use~\eqref{eq:mod_rec_no_ui'} starting with
\begin{multline}
 \left[ J^{(b)}_j, M^{(1,3)}_i  \right] = -b^2 (i-1)(i-2)(i-3) \delta_{i+j,3} - b j(j-1) J^{(b)}_{i+j-3}\delta_{j \leq 1} \\ 
+ b(i-1)j J^{(b)}_{i+j-3} \left( \delta_{i+j \geq 3 } + \delta_{j \leq 0 } \right) + (1+b)j M^{(1,2)}_{i+j-1}\delta_{i+j\geq 2} \\
+ \sum\limits_{n\geq j} (1+b) j J^{(b)}_{i+j-n-2} M^{(1,1)}_{n}\left( \delta_{n \geq1} +  \delta_{j\leq 0} \right).
\end{multline}
For $j\geq 1$, this leads to
\begin{equation}
[p_j^*,  M^{(1,3)}_i] = j\left(\sum_{l\geq 0} + \sum_{l\geq j-1}\right) J^{(b)}_{i+j-3-l} J^{(b)}_l + \frac{b}{1+b}j(2i+j-3) J^{(b)}_{i+j-3} \delta_{i+j\geq 3}
\end{equation}
and therefore (recall that $p^*_n = 0$ if $n\leq 0$)
\begin{multline} \label{Commutatorp*M13}
[p_i^*,  M^{(1,3)}_j] - [p_j^*,  M^{(1,3)}_i] = b(i-j)(i+j-3) p^*_{i+j-3}\\
+ \frac{1}{1+b}\left(2(i-j) \sum_{l\geq M-1} -\sgn(i-j) \mu \sum_{l=\mu-1}^{M-2} + (i-j)\sum_{l=0}^{M-2}\right)J^{(b)}_{i+j-3-l} J^{(b)}_l
\end{multline}
and (recall that $M^{(1,1)}_n = 0$ if $n\leq 0$)
\begin{multline} \label{CommutatorM11M13}
[M^{(1,1)}_i,  M^{(1,3)}_j] - [M^{(1,1)}_j,  M^{(1,3)}_i] = b(i-j)(i+j-3) M^{(1,1)}_{i+j-3}\\
\left(2(i-j) \sum_{l\geq M-1} -\sgn(i-j) (\mu-1) \sum_{l=\mu-1}^{M-2} + (i-j)\sum_{l=1}^{M-2}\right)J^{(b)}_{i+j-3-l} M^{(1,1)}_l.
\end{multline}
The RHS of those equations should be respectively equal to 
\begin{multline} \label{DTilde3p*}
\sum_{l\geq 1} \tilde{D}_{ij, l}(3) p_l^* = b(i-j)(i+j-3) p_{i+j-3}^*\\
+ \left(2(i-j)\sum_{l\geq M-1} + (i-j) \sum_{l=M-1}^{i+j-3} + \sgn(i-j) \sum_{l=\mu-1}^{M-2} (2l-3\mu+3)\right) J^{(b)}_{i+j-3-l} p_l^*
\end{multline}
and
\begin{multline}
\sum_{l\geq 1} \tilde{D}_{ij, l}(3) M^{(1,1)}_l + \tilde{D}_{ij, l}(1) M^{(1,3)}_l = b(i-j)(i+j-3) M^{(1,1)}_{i+j-3}\\
+ \left(2(i-j)\sum_{l\geq M-1} + (i-j) \sum_{l=M-1}^{i+j-3} + \sgn(i-j) \sum_{l=\mu-1}^{M-2} (2l-3\mu+3)\right) J^{(b)}_{i+j-3-l} M^{(1,1)}_l
\end{multline}
This can be proven exactly as we did for 3-constellations in Section \ref{sec:CommutatorA1A3}, and therefore we do not repeat it. There is however a subtlety in equating \eqref{Commutatorp*M13} with \eqref{DTilde3p*}, because the former is written with $J^{(b)}_l$, where $l=0$ is not vanishing, instead of $p_l^*$ for \eqref{DTilde3p*}. To deal with this, one can separate three cases: $i)$ $\mu>1$, $ii)$ $\mu=1$ and $M>2$, $iii)$ $\mu=1, M=2$ and prove the equality in all of them.

\subsection{Computation of $[M^{(1,2)}_i, M^{(1,3)}_j]$}
By writing again $M^{(1,3)}_j = \sum_{l\geq 1} J^{(b)}_{j-l-1} M^{(1,2)}_l + b(j-1) M^{(1,2)}_{j-1}$, one finds
\begin{multline}
[M^{(1,2)}_i, M^{(1,3)}_j] =  b((i-1)(i-2) + (j-1)(i-j+1)) M^{(1,2)}_{i+j-3}\\
+\left(\sum_{n\geq j-1} (n-j+1) + \sum_{n=1}^{i+j-3} (n-j+1) + \sum_{n\geq i-1} (2i-n-2)\right) J^{(b)}_{i+j-3-n} M^{(1,2)}_{n}.
\end{multline}
It is then elementary to write
\begin{multline}
[M^{(1,2)}_i, M^{(1,3)}_j] - [M^{(1,2)}_j, M^{(1,3)}_i] = (i-j) M^{(1,3)}_{i+j-2} + b(i-j)(i+j-3) M^{(1,2)}_{i+j-3} \\
+ \left(2(i-j)\sum_{n\geq M-1} + (i-j) \sum_{n=M-1}^{i+j-3} + \operatorname{sgn}(i-j)\sum_{n=\mu-1}^{M-2}(2n-3\mu+3)\right) J^{(b)}_{i+j-3-n} M^{(1,1)}_n,
\end{multline}
i.e. $[M^{(1,2)}_i, M^{(1,3)}_j] - [M^{(1,2)}_j, M^{(1,3)}_i] = \sum_{n\geq 1} \tilde{D}_{ij, n}(2) M^{(1,3)}_n + \tilde{D}_{ij, n}(3) M^{(1,2)}_n$ as desired.

\subsection{Computation of $\left[ M^{(1,3)}_i, M^{(1,3)}_j\right]$}
It follows the same steps as $[A_i(3), A_j(3)]$ for 3-constellations in Section \ref{sec:A3A3}. Since it is a key result of this article, we provide here some detail. In order to emphasize the similarity with $[A_i(3), A_j(3)]$, we decompose the calculation of $[M^{(1,3)}_i, M^{(1,3)}_j]$ in the same way and to avoid new notation, {\bf we reallocate the notation of Section \ref{sec:A3A3}} to the analogous quantities they represent in our decomposition of $[M^{(1,3)}_i, M^{(1,3)}_j]$. The difference in the notation reallocation is that the summands are now of the form $JJM, JM, bJM, \dotsc$

The commutator $[ M^{(1,3)}_i, M^{(1,3)}_j]$ expands as
{\small 
\begin{multline}
[ M^{(1,3)}_i, M^{(1,3)}_j] = \underbrace{\sum\limits_{\substack{\ell \geq 1 }} \sum\limits_{\substack{n \geq 1 }} \left[ J^{(b)}_{i-1-\ell}M^{(1,2)}_\ell,J^{(b)}_{j-1-l}M^{(1,2)}_n\right]}_{Q_{ij}} + \underbrace{b^2(j-1)(i-1) \left[M^{(1,2)}_{i-1},M^{(1,2)}_{j-1} \right]}_{H_{ij}} \nonumber\\
+ \underbrace{b(j-1)\sum\limits_{\substack{\ell \geq 1 }} \left[ J^{(b)}_{i-1-\ell}M^{(1,2)}_\ell,M^{(1,2)}_{j-1}\right] }_{G_{ij}} +\underbrace{b(i-1)\sum\sum\limits_{\substack{n \geq 1 }} \left[M^{(1,2)}_{i-1},J^{(b)}_{j-1-n}M^{(1,2)}_n\right] }_{-G_{ji}} .
\end{multline}}%
We decompose $Q_{ij}$ by expanding the commutator
\begin{multline}
Q_{ij} = \underbrace{\sum\limits_{\substack{\ell \geq 1 }} \sum\limits_{\substack{n \geq 1 }} J^{(b)}_{i-1-\ell} [M^{(1,2)}_\ell, J^{(b)}_{j-1-n}] M^{(1,2)}_n}_{P_{ij}} + \underbrace{\sum\limits_{\substack{\ell \geq 1 }} \sum\limits_{\substack{n \geq 1 }} J^{(b)}_{j-1-n} [J^{(b)}_{i-1-\ell} , M^{(1,2)}_n] A^2_\ell}_{-P_{ji}} \nonumber\\
+ \underbrace{\sum\limits_{\substack{\ell \geq 1 }} \sum\limits_{\substack{n \geq 1 }} [J^{(b)}_{i-1-\ell} , J^{(b)}_{j-1-n}] M^{(1,2)}_n M^{(1,2)}_\ell}_{B_{ij}} + \underbrace{\sum\limits_{\substack{\ell \geq 1 }} \sum\limits_{\substack{n \geq 1 }} J^{(b)}_{i-1-\ell} J^{(b)}_{j-1-n} [ M^{(1,2)}_\ell, M^{(1,2)}_n]}_{C_{ij}}. 
\end{multline}
We now rewrite each of these contributions individually before assembling them to show that the commutator closes.

\subsubsection{Term $B_{ij}$} This term writes
\begin{equation}
B_{ij} = (1+b) \sum\limits_{\ell = 1}^{i+j-3} (i-1-\ell) M^{(1,2)}_{i+j-2-\ell} M^{(1,2)}_\ell .
\end{equation}
Now commuting the two factors $M^{(1,2)}$, relabelling the sum using $\ell\to i+j-2-\ell$ and taking half-sums of the two expressions gives
\begin{multline}
B_{ij} = \underbrace{\frac{1+b}{12} (i+j-2)(i+j-3)(i+j-4) M^{(1,2)}_{i+j-4}}_{B^{M}_{ij}} \\+ \underbrace{\frac{1+b}{2}(i-j) \sum\limits_{\ell=1}^{i+j-3} M^{(1,2)}_{i+j-2-\ell} M^{(1,2)}_\ell}_{B^{MM}_{ij}}.
\end{multline}

\subsubsection{Term $C_{ij}$} It is given by
\begin{equation}
    C_{ij} = \sum\limits_{\substack {\ell \geq 1 }} \sum\limits_{\substack {n \geq 1 }} (\ell-n) J^{(b)}_{i-1-\ell}  J^{(b)}_{j-1-n} M^{(1,2)}_{\ell+n-2} = \sum\limits_{\substack{ n \geq 0 }} \sum\limits_{\substack{\ell = j-1 }}^{n+j-1}  (2\ell-2j-n+2) J^{(b)}_{i+j-3-\ell} J^{(b)}_{\ell-1-n} M^{(1,2)}_n.
\end{equation}
Following the exact same steps as in Section \ref{sec:C} we arrive at
\begin{equation}
C_{ij} = C^{JJM,1}_{ij} + C^{JJM,2}_{ij} + C_{ij}^{M,1} + C^{M,2}_{ij} + \tilde{C}^{JJM}_{ij}
\end{equation}
with
\begin{equation}
\begin{aligned}
C^{JJM,1}_{ij} &= (i-j) \sum_{ \substack{\ell \geq M-\mu}} \sum_{\substack{n=M-1 }}^{\ell+\mu-1} J^{(b)}_{i+j-3-n} J^{(b)}_{n-1-\ell} A^2_\ell\\
C^{JJM,2}_{ij} &= (i-j)\sum_{\ell =0}^{M-\mu-1} \sum_{\substack{ n=\mu+\ell}}^{M-2} J^{(b)}_{i+j-3-n} J^{(b)}_{n-1-\ell}A^2_\ell\\
C_{ij}^{M,1} &= - B^{M}_{ij}\\
C^{M,2}_{ij} &= (1+b)\frac{\operatorname{sgn}(j-i)}{2}\sum_{n=\mu-1}^{M-2} (2n+6-3\mu-M)(i+j-3-n) M^{(1,2)}_{i+j-4}\\
\tilde{C}^{JJM}_{ij} &= \operatorname{sgn}(j-i) \sum_{\substack{\ell \geq 1 }} \sum_{n=\mu-1}^{M-2} (2\mu-2n+\ell-2) J^{(b)}_{i+j-3-n} J^{(b)}_{n-1-\ell} M^{(1,2)}_\ell
\end{aligned}
\end{equation}

\subsubsection{Term $P_{ij}$}
This term expands as 
\begin{multline}
    P_{ij} = \underbrace{ b \sum\limits_{n \geq j-1} (n-j)(n-j+1) J^{(b)}_{i+j-3-n} M^{(1,2)}_{n-1}}_{F_{ij}}
    + \underbrace{\sum\limits_{n,\ell \geq j-1} (\ell-j+1) J^{(b)}_{i+j-3-n} J^{(b)}_{n-1-\ell} A^2_{\ell}}_{D_{ij}} \nonumber\\
    + \underbrace{\sum\limits_{n \geq j-1} \sum\limits_{ \substack{ \ell =1 }}^{n-1} (\ell-j+1) J^{(b)}_{i+j-3-n} J^{(b)}_{n-1-\ell} A^2_{\ell}}_{E_{ij}}. \nonumber
\end{multline}
Following the same steps as in Section \ref{sec:P}, we arrive at
\begin{multline}
D_{ij}-D_{ji}
= (i-j) \sum\limits_{\substack{n \geq M-1 }} J^{(b)}_{i+j-3-n} M^{(1,3)}_n \\
+ D^{MM,1}_{ij} + D^{MM,2}_{ij} + D^M_{ij} + D^{JJM}_{ij} + \tilde{D}^{JJM,1}_{ij} + \tilde{D}^{JJM,2}_{ij} + D^{bJM,1}_{ij} + D^{bJM,2}_{ij}
\end{multline}
with
\begin{equation}
\begin{aligned}
D^{MM,1}_{ij} &= -(1+b)\frac{(i-j)}{2}\sum\limits_{\ell = 1}^{M+\mu-3}  M^{(1,2)}_{i+j-2-\ell} M^{(1,2)}_\ell\\
D^{MM,2}_{ij} &= -(1+b)\frac{(i-j)}{2}\sum\limits_{\ell = \mu}^{M-2}  M^{(1,2)}_{i+j-2-\ell} M^{(1,2)}_\ell\\
D^M_{ij} &= (1+b)\frac{(i-j)}{2}\sum\limits_{\substack{\ell =1}}^{M-2} (i+j-2-2\ell) M^{(1,2)}_{i+j-4}\\
D^{JJM}_{ij} &= (i-j)\sum\limits_{\substack{\ell =1}}^{M-2} \sum\limits_{\substack{ n=\ell+1 }}^{M-2}  J^{(b)}_{i+j-3-n} J^{(b)}_{n-1-\ell} M^{(1,2)}_\ell\\
\tilde{D}^{JJM,1}_{ij} &= \sgn(j-i)\sum\limits_{\substack{n = \mu-1}}^{M-2} \sum\limits_{\substack{\ell \geq \mu-1 }} (\mu-1-\ell) J^{(b)}_{i+j-3-n} J^{(b)}_{n-1-\ell} M^{(1,2)}_\ell\\
\tilde{D}^{JJM,2}_{ij} &= \sgn(j-i) \sum\limits_{\substack{n \geq M-1 }} \sum\limits_{\substack{\ell = \mu-1 }}^{M-2} (\mu-1-\ell) J^{(b)}_{i+j-3-n} J^{(b)}_{n-1-\ell} M^{(1,2)}_\ell\\
D^{bJM,1}_{ij} &= - b(i-j) \sum\limits_{\substack{n \geq M-1 }} (n-1) J^{(b)}_{i+j-3-n} M^{(1,2)}_{n-1}\\
D^{bJM,2}_{ij} &= b(i-j) \sum\limits_{\substack{\ell =1}}^{M-2} (i+j-3-\ell) J^{(b)}_{i+j-4-\ell-n} M^{(1,2)}_\ell
\end{aligned}
\end{equation}
for $D_{ij}-D_{ji}$. As for $E_{ij}-E_{ji}$,
\begin{equation}
E_{ij} - E_{ji} = (i-j)\sum\limits_{\substack{ n \geq M-1 }} J^{(b)}_{i+j-3-n} M^{(1,3)}_{n} + E^{JJM}_{ij} + \tilde{E}^{JJM,1}_{ij} + \tilde{E}^{JJM,2}_{ij} + E^{bJM}_{ij},
\end{equation}
with
\begin{equation}
\begin{aligned}
E^{JJM}_{ij} &= -(i-j) \sum\limits_{\substack{ n \geq M-1 }} \sum\limits_{ \ell \geq n-1} J^{(b)}_{i+j-3-n} J^{(b)}_{n-1-\ell} M^{(1,2)}_\ell\\
\tilde{E}^{JJM,1}_{ij} &= \sgn(j-i) \sum\limits_{\ell=1}^{\mu-2} \sum\limits_{\substack{n = \mu-1 }}^{M-2} (\mu-1-\ell) J^{(b)}_{i+j-3-n} J^{(b)}_{n-1-\ell} M^{(1,2)}_\ell\\
\tilde{E}^{JJM,2}_{ij} &=  \sgn(j-i) \sum\limits_{\ell=\mu-1}^{M-3} \sum\limits_{\substack{n=\ell+1}}^{M-2} (\mu-1-\ell)   J^{(b)}_{i+j-3-n} J^{(b)}_{n-1-\ell} M^{(1,2)}_\ell\\
E^{bJM}_{ij} &= -b(i-j)\sum\limits_{\substack{ n \geq M-1 }} (n-1) J^{(b)}_{i+j-3-n} M^{(1,2)}_{n-1}
\end{aligned}
\end{equation}
Finally,
\begin{multline}
F_{ij}-F_{ji} = \underbrace{b(i-j)\sum\limits_{n \geq M-1} (2n-i-j+1)J^{(b)}_{i+j-3-n} M^{(1,2)}_{n-1}}_{F^{bJM}_{ij}} \\
+ \underbrace{b\sgn(i-j)\sum\limits_{n = \mu-1}^{M-2} (\mu-n)(\mu-1-n)J^{(b)}_{i+j-3-n} M^{(1,2)}_{n-1}}_{\tilde{F}^{bJM}_{ij}}.
\end{multline}

\subsubsection{Terms $G_{ij}$ and $H_{ij}$} The term $H_{ij} = b^2(j-1)(i-1) [M^{(1,2)}_{i-1},M^{(1,2)}_{j-1}]$ is direct to compute and gives
\begin{equation}
H_{ij} = b^2(i-1)(j-1)(i-j) M^{(1,2)}_{i+j-4}.
\end{equation}

The term $G_{ij} = b(j-1)\sum\limits_{\substack{\ell \geq 1 }} [ J^{(b)}_{i-1-\ell}M^{(1,2)}_\ell,M^{(1,2)}_{j-1}]$ gives
\begin{multline}
    G_{ij}-G_{ji} = \underbrace{b^2\left[(i-1)(i-2)(i-3)-(j-1)(j-2)(j-3)\right] M^{(1,2)}_{i+j-4}}_{G^{M}_{ij}} \\
    + \underbrace{b(i-j)\sum\limits_{n=1}^{i+j-4} n J^{(b)}_{i+j-4-n} M^{(1,2)}_{n}}_{G^{bJM,1}_{ij}}
    + \underbrace{2b(i-j)\sum\limits_{n\geq M-2} (i+j-3)J^{(b)}_{i+j-4-n}M^{(1,2)}_n}_{G^{bJM,2}_{ij}} 
    + \tilde{G}^{bJM}_{ij}
\end{multline}
with 
\begin{multline}
\tilde{G}^{bJM}_{ij} = b\sgn(i-j)\sum\limits_{n=\mu-2}^{M-3} \left((\mu-1)(n-2\mu+4)-(M-1)(\mu-1-n)\right)J^{(b)}_{M+\mu-4-n}M^{(1,2)}_n \\
+ b\sgn(i-j)\left((M-1) J^{(b)}_{M-2} M^{(1,2)}_{\mu-2} - (\mu-1) J^{(b)}_{\mu-2} M^{(1,2)}_{M-2}\right)
\end{multline}

%

\subsubsection{Repackaging contributions}
We start by adding together $\tilde{C}_{ij}^{JJM}$, $\tilde{D}_{ij}^{JJM,1}$ and $\tilde{E}_{ij}^{JJM,1}$ 
\begin{equation}
\begin{aligned}
\tilde{C}_{ij}^{JJM} + \tilde{D}_{ij}^{JJM,1} + \tilde{E}_{ij}^{JJM,1} &= \sgn(j-i) \sum\limits_{n=\mu-1}^{M-2} \sum\limits_{\ell \geq 1} (3\mu-3-2n) J^{(b)}_{i+j-3-n} J^{(b)}_{n-1-\ell} M^{(1,2)}_\ell \\
&\begin{multlined} = \sgn(j-i) \sum\limits_{n=\mu-1}^{M-2}(3\mu-3-2n) J^{(b)}_{i+j-3-n} M^{(1,3)}_n \\
 \underbrace{-b\sgn(j-i) \sum\limits_{n=\mu-1}^{M-2}(3\mu-3-2n) (n-1) J^{(b)}_{i+j-3-n} M^{(1,2)}_{n-1}}_{{\widetilde{(CE)}}^{bJM}_{ij}}.\end{multlined}
\end{aligned}
\end{equation}

Now we add together the contributions $\tilde{D}_{ij}^{JJM,2}$ and $\tilde{E}_{ij}^{JJM,2}$. They have similar range except for the $\ell = M-2$ contribution of  $\tilde{D}_{ij}^{JJM,2}$ which is
\begin{equation}
    \lambda_{ij} \coloneqq \sgn(j-i)\sum\limits_{n\geq M} (\mu-M+1)J^{(b)}_{M+\mu-3-n} J^{(b)}_{n+1-M} M^{(1,2)}_{M-2}.
\end{equation}
Now the sum over $n$ in $\tilde{D}_{ij}^{JJM,2}+\tilde{E}_{ij}^{JJM,2}$ can be encoded into some $M^{(1,2)}$ as follows
\begin{align}
\tilde{D}_{ij}^{JJM,2}+\tilde{E}_{ij}^{JJM,2} &=  \lambda_{ij} + \sgn(j-i)(1+b)\sum\limits_{\ell=\mu}^{M-3} (\mu-1-\ell) M^{(1,2)}_{M+\mu-2-\ell} M^{(1,2)}_{\ell} \\
&-b(1+b)\sgn(j-i)\sum\limits_{\ell=\mu-1}^{M-3} (\mu-1-\ell)(M+\mu-3-\ell) M^{(1,1)}_{M+\mu-3-\ell}M^{(1,2)}_{\ell}. \nonumber
\end{align}
Commuting the two $M^{(1,2)}$s and relabelling the sum as $\ell \rightarrow M+\mu-2-\ell$ then taking the half-sum gives
\begin{multline}
\sgn(j-i) \sum\limits_{\ell=\mu}^{M-3} (\mu-1-\ell) M^{(1,2)}_{M+\mu-2-\ell} M^{(1,2)}_{\ell} = \frac{i-j}{2} \sum_{\ell=\mu+1}^{M-3} M^{(1,2)}_{i+j-2-\ell} M^{(1,2)}_\ell\\
+ \frac{\sgn(j-i)}{2} \left(\sum_{\ell=\mu+1}^{M-2} (\ell+1-M)(2\ell +2-i-j) M^{(1,2)}_{i+j-4} - M^{(1,2)}_\mu M^{(1,2)}_{M-2} -  M^{(1,2)}_{M-2} M^{(1,2)}_\mu\right)
\end{multline}
Similarly, the sum over $n$ in $\lambda_{ij}$ can also be written in terms of some $M^{(1,2)}$,
\begin{multline}
\lambda_{ij} = \frac{1}{2} (i-j+\sgn(j-i)) (M^{(1,2)}_\mu M^{(1,2)}_{M-2} + M^{(1,2)}_{M-2} M^{(1,2)}_\mu + (\mu-M+2) M^{(1,2)}_{i+j-4}) \\
-b(1+b)\sgn(j-i)(M-\mu-1)(\mu-1) M^{(1,1)}_{\mu-1} M^{(1,2)}_{M-2}
\end{multline}
%
Altogether, one gets
\begin{multline}
\tilde{D}_{ij}^{JJM,2}+\tilde{E}_{ij}^{JJM,2} = \underbrace{\frac{1+b}{2}(i-j) \sum\limits_{\ell=\mu}^{M-2} M^{(1,2)}_{i+j-2-\ell} M^{(1,2)}_{\ell}}_{(DE)_{ij}^{MM}}\\
+ \underbrace{(1+b)\frac{\sgn(j-i)}{2}\sum\limits_{\ell=\mu}^{M-2} (\ell-M+1)(2\ell-i-j+2) M^{(1,2)}_{i+j-4}}_{(DE)^M_{ij}}\\
+ \underbrace{ b(1+b)\sgn(i-j)\sum\limits_{\ell=\mu-1}^{M-2} (\mu-1-\ell)(i+j-3-\ell) M^{(1,1)}_{i+j-3-\ell}M^{(1,2)}_{\ell}}_{\widetilde{(DE)}_{ij}^{bJM}}
\end{multline}

Next, we add together the terms $C^{JJM,1}_{ij}$ and $E^{JJM}_{ij}$. The same treatment as for $C^{JJA,1}_{ij} + E^{JJA}_{ij}$ in Section \ref{sec:Repackage} leads to
\begin{equation}
C^{JJM,1}_{ij} + E^{JJM}_{ij} = (i-j)\sum\limits_{n=M-1}^{i+j-3}  J^{(b)}_{i+j-3-n} M^{(1,3)}_{n} + (CE)^{JJM,1}_{ij} + (CE)^{JJM,2}_{ij} + (CE)^M_{ij} + (CE)^{bJM}_{ij}
\end{equation}
with
\begin{equation}
\begin{aligned}
(CE)^{JJM,1}_{ij} &= (i-j)\sum\limits_{n=M-1}^{i+j-3} \sum\limits_{\ell = n-\mu+1}^{M-2}  J^{(b)}_{i+j-3-n} J^{(b)}_{n-1-\ell} M^{(1,2)}_{\ell}\\
(CE)^{JJM,2}_{ij} &= -(i-j) \sum\limits_{n=M-1}^{i+j-3} \sum\limits_{\ell = 1}^{M-2} J^{(b)}_{i+j-3-n} J^{(b)}_{n-1-\ell} M^{(1,2)}_{\ell}\\
(CE)^M_{ij} &= -(1+b)(i-j)\frac{(\mu-1)(\mu-2)}{2} M^{(1,2)}_{i+j-4}\\
(CE)^{bJM}_{ij} &= -b(i-j)\sum\limits_{n=M-1}^{i+j-3} (n-1)J^{(b)}_{i+j-3-n} M^{(1,2)}_{n-1}
\end{aligned}
\end{equation}

Finally, we repackage together
\begin{align}
G^{M}_{ij} + H_{ij} = b^2(i-j)(i+j-3)(i+j-4) M^{(1,2)}_{i+j-4}.
\end{align}
Moreover,
\begin{multline}
D^{bJM,1}_{ij} + E^{bJM}_{ij} + F^{bJM}_{ij} + G^{bJM,2}_{ij} = b(i-j)(i+j-3) \sum_{n\geq M-1} J^{(b)}_{i+j-3-n} M^{(1,2)}_{n-1} 
\end{multline}
and
\begin{multline}
(CE)^{bJM}_{ij} + G^{bJM,1}_{ij} + D^{bJM,2}_{ij} = b(i-j)(i+j-3) \sum_{n=1}^{M-2} J^{(b)}_{i+j-3-n} M^{(1,2)}_{n-1} \\ + \underbrace{b(i-j)(\mu-1) J^{(b)}_{\mu-2} M^{(1,2)}_{M-2}}_{\theta_{ij}}. 
\end{multline}
Hence, adding together the last three equations gives a contribution
\begin{multline}
(D^{bJM,1}_{ij} + E^{bJM}_{ij} + F^{bJM}_{ij} + G^{bJM,2}_{ij}) + ((CE)^{bJM}_{ij} + G^{bJM,1}_{ij} + D^{bJM,2}_{ij}) + (H_{ij} + G^{M}_{ij}) \\
= b(i-j)(i+j-3) M^{(1,3)}_{i+j-3} + \theta_{ij}.
\end{multline}

\subsubsection{Cancellations}
We notice the following cancellations,
\begin{equation}
\begin{aligned}
C^{M,1}_{ij} + B^{M}_{ij} &= 0, \\
D^{MM,1}_{ij} + B^{MM}_{ij} &= 0, \\
D^{MM,2}_{ij} + (DE)^{MM}_{ij} &=0, \\
C^{M,2}_{ij} + D^M_{ij} + (DE)^M_{ij} + (CE)^M_{ij}  &= 0, \\
C^{JJM,2}_{ij} + D^{JJM,1}_{ij} + (CE)^{JJM,1}_{ij} + (CE)^{JJM,2}_{ij} &= 0, \\
\widetilde{(CE)}^{bJM}_{ij} + \widetilde{(DE)}^{bJM}_{ij} + \tilde{G}^{bJM}_{ij} + \tilde{F}^{bJM}_{ij} + \theta_{ij} &= 0.
\end{aligned}
\end{equation}

\subsubsection{Final expression} Bringing together the remaining terms, we get the final expression for the commutator
\begin{multline}
\label{eq:comm_triv_bip}
[M^{(1,3)}_i, M^{(1,3)}_j] = 2(i-j) \sum\limits_{n \geq M-1} J^{(b)}_{i+j-3-n} M^{(1,3)}_n + (i-j) \sum\limits_{n=M-1}^{i+j-3} J^{(b)}_{i+j-3-n} M^{(1,3)}_n \\
+ \sgn(j-i)\sum\limits_{n=\mu-1}^{M-2} (3\mu-2n-3) J^{(b)}_{i+j-3-n} M^{(1,3)}_n 
+b(i-j)(i+j-3) M^{(1,3)}_{i+j-3} 
\end{multline}
which concludes the computation and the proof of Proposition~\ref{thm:SimplifiedCommutators2}.

\bibliography{bDeformedConstraints}
\bibliographystyle{alpha}

\end{document}